\documentclass[leqno,A4]{aadmbook}
\usepackage{amsthm}
\usepackage{amsfonts}
\usepackage{amsmath}
\usepackage{graphicx}
\usepackage{float} %package for figure and table placement 
\usepackage{amsmath}

\newcommand{\AAA}{{\mathbb A}}
\newcommand{\BBB}{{\mathbb B}}
\newcommand{\SSS}{{\mathbb S}}
\newcommand{\RR}{{\mathbb R}}

\newcommand{\NN}{{\mathbb N}}

\newtheorem{statement}{Statement}
\newtheorem{theorem}{Theorem}

\newtheorem{lemma}{Lemma}

\newtheorem{proposition}{Proposition}
\newtheorem{remark}{Remark}

\begin{document}

% Your \newcommands below (if there are any):

\oddsidemargin 16.5mm
\evensidemargin 16.5mm

\thispagestyle{plain}

\begin{center}
{\large \sc  Applicable Analysis and Discrete Mathematics}

{\small available online at  http:/$\!$/pefmath.etf.rs }
\end{center}

\noindent{\small{\sc  Appl. Anal. Discrete Math.\ }{\bf 18} (2024), 244--288.}\\
\noindent{\scriptsize https://doi.org/10.2298/AADM240308012B}     %%%%%%%%%%%%%%%%%%%%%%%%%%%%%%%%%%%%%%%%%%%%

\vspace{5cc}
\begin{center}

{\large\bf  THE BEST POSSIBLE CONSTANTS APPROACH \\FOR WILKER-CUSA-HUYGENS INEQUALITIES \\ [-1.1ex] VIA STRATIFICATION
\rule{0mm}{6mm}\renewcommand{\thefootnote}{}%Enter at least one, but not more than 3 MSCs.
% First entered MSC will be a primary one, others (at most 2) will be secondary.
\footnotetext{\scriptsize ${}^{\ast}$Corresponding author. Milo\v s Mi\' covi\' c}
\footnotetext{

\rule{2.4mm}{0mm}Keywords and Phrases. Wilker-Cusa-Huygens inequalities, Stratified families of functions, A mini-
\rule{2.4mm}{0mm}\indent max approximant, Automated proving of MTP inequalities, SimTheP.}}

\vspace{1cc}
{\large\it Bojan Banjac, Branko Male\v sevi\' c, Milo\v s Mi\' covi\' c$\,{}^{\ast}$, \\ Bojana Mihailovi\' c and Milica Savatovi\' c}

\vspace{1cc}
\parbox{24cc}{{\small

In this paper, we generalize Cristinel Mortici's results on Wilker-Cusa-Huygens inequalities using stratified families of functions and SimTheP -- a system for automated proving of MTP inequalities.

}}
\end{center}

\vspace{1.5cc}
\begin{section}
{INTRODUCTION}
\end{section}

The basis of this research is well-known C. Mortici's paper \cite{Mortici_2011}
in which the following theorems were proved:

\begin{theorem}
For every $x \in (0,\pi/2)$, we have$:$
$$
-\dfrac{1}{15}x^4
<
\cos x - \left(\dfrac{\sin x}{x}\right)^{\!3}
<
-\dfrac{1}{15}x^4 + \dfrac{23}{1890}x^6 \,.
$$
\end{theorem}

\begin{theorem}
For every $x \in (0,\pi/2)$, we have$:$
$$
-\dfrac{1}{180}x^4
<
\dfrac{\sin x}{x} - \dfrac{\cos x + 2}{3}
<
-\dfrac{1}{180}x^4 + \dfrac{1}{3780} x^6 \,.
$$
\end{theorem}

\begin{theorem}
For every $x \in (0,\pi/2)$, we have$:$
$$
3 + \left(\dfrac{3}{20}x^4 - \dfrac{3}{140}x^6 \right)\dfrac{1}{\cos x}
<
2\dfrac{\sin x}{x} + \dfrac{\tan x}{x}
<
3 + \dfrac{3}{20 \cos x}x^4 \,.
$$
\end{theorem}

\begin{theorem}
For every $x \in (0,\pi/2)$, we have$:$
$$
2 + \left( \dfrac{8}{45}x^4 - \dfrac{4}{105}x^6 \right) \dfrac{1}{\cos x}
<
\left(\dfrac{\sin x}{x}\right)^2 + \dfrac{\tan x}{x}
<
2 + \dfrac{8}{45 \cos x}x^4 \,.
$$
\end{theorem}

\begin{theorem}
For every $x \in (0,\pi/2)$, we have$:$
$$
\left(\dfrac{x}{\sin x}\right)^2 + \dfrac{x}{\tan x} > 2 + \dfrac{2}{45}x^4 \,.
$$
\end{theorem}

\begin{theorem} 
For every $x \in (0,\pi/2)$, we have$:$
$$
3\dfrac{x}{\sin x} + \cos x > 4 + \dfrac{1}{10}x^4 + \dfrac{1}{210}x^6 \,.
$$
\end{theorem}

Using a method based on the stratified families of functions described in the paper \cite{Malesevic_Mihailovic_2021}, we show that it is possible to enhance inequalities presented in Theorems 1-6 up to the level of the best possible constants. Also, we show that from some of the introduced families of functions, it is possible to single out corresponding minimax approximations.

In Section 2, preliminaries on stratified families of functions and a method for proving MTP inequalities are presented. This method for proving MTP inequalities forms the basis of SimTheP -- automated theorem prover for MTP inequalities. The main results of the paper are provided in Section 3. Section 4 presents the conclusion. At the end of the paper, an Appendix is included, which contains proofs obtained using SimTheP.

\vspace{1.5cc}
\begin{section}
{PRELIMINARIES}
\end{section}

In this section, we present assertions from the paper \cite{Malesevic_Mihailovic_2021} and a method from the paper \cite{Malesevic_Makragic_2016} on the basis of which our main results will be obtained.

\bigskip\noindent
{\bf 2.1. STRATIFIED FAMILIES OF FUNCTIONS}

\medskip

Let
$$
\varphi_p(x):(a,b)\longrightarrow \mathbb{R}
$$
be a family of functions with a variable $x \in (a,b)$ and a parameter $p\in \mathbb{R}^{+}$.

A family of functions $\varphi_p(x)$ is
{\em increasingly stratified} if
$$(\forall p_1,p_2 \in \RR^{+}) \,\,\, p_1 \!<\! p_2 \Longleftrightarrow
\varphi_{p_1}(x) \!<\! \varphi_{p_2}(x)$$
holds for any $x \!\in\! (a,b)$ and, conversely, it is
{\em decreasingly stratified} if
$$(\forall p_1,p_2 \!\in\! \RR^{+}) \,\,\, p_1 \!<\! p_2\Longleftrightarrow
\varphi_{p_1}(x) \!>\! \varphi_{p_2}(x)$$
holds for any $x \!\in\! (a,b)$.

\break

In this paper, we call $\sup\limits_{x \in \left(a, b\right)} \left| \varphi_{p}(x) \right|$ {\em an error} and denote it by $d^{(p)}$.

In \cite{Malesevic_Mihailovic_2021}, the conditions for the existence of a unique value $p_0$ of the parameter $p \in \mathbb{R}^{+}$, for which an infimum of the error is attained, are considered. Such infimum is denoted by
$$
d_0 = \inf\limits_{p \in \RR^{+}} \sup\limits_{x \in (a,b)}{\left| \varphi_p(x) \right|}.
$$
For such a value $p_0$, the function $\varphi_{p_0}(x)$ is called {\em the minimax approximant} on $(a,b)$.

\begin{theorem} {\rm (Theorem 1 \cite{Malesevic_Mihailovic_2021})}
\label{Existence_and_Uniqueness}
Let $\varphi_p(x)$ be a family of functions that are continuous with respect to $x \!\in\! (a,b)$
for each $p \!\in\! \RR^{+}$ and increasingly stratified for $p \!\in\! \RR^{+}$, and let $c,d$ be in $\RR^{+}$
such that $c<d$. If$\,:$
\begin{itemize}
\item[{\rm (a)}] $\varphi_c(x) \!<\! 0$ and $\varphi_d(x) \!>\! 0$ for all $x \!\in\! (a,b)$,
    and at the endpoints $\varphi_c(a+) = \varphi_c(b-) = \varphi_d(a+) = 0$ and $\varphi_d(b-) \!\in\! \RR^{+}$ hold;
\item[{\rm (b)}] the functions $\varphi_p(x)$ are continuous
    with respect to $p \!\in\! (c,d)$ for each $x \!\in\! (a,b)$
    and $\varphi_p(b-)$ is continuous with respect to $p \!\in\! (c,d)$ too;
\item[{\rm (c)}] for all $p \!\in\! (c,d)$, there exists a right neighbourhood of the point $a$ in which
    \mbox{$\varphi_p(x) \!<\! 0$} holds;
\item[{\rm (d)}] for all $p \!\in\! (c,d)$  the function $\varphi_p(x)$ has exactly one extremum $t^{(p)}$ on $(a,b)$, which is minimum;
\end{itemize}
then there exists exactly one solution $p_0$, for $p \!\in\! \RR^{+}$, of the following equation
$$
|\varphi_p(t^{(p)})|
=
\varphi_p(b-)
$$
and for $d_0=|\varphi_{p_0}(t^{(p_0)})|=\varphi_{p_0}(b-)$ we have
$$
d_0
\, =
\inf\limits_{\mbox{\scriptsize $p \!\in\! \RR^{+}$}}
\sup\limits_{\mbox{\scriptsize $x \!\in\! (a,b)$}} |\varphi_p(x)|\,.
$$
\end{theorem}

The analogous theorem can be stated for decreasingly stratified families of functions.

\medskip
\noindent
{\bf Theorem 7'} {\rm (Theorem  1' \cite{Malesevic_Mihailovic_2021})} \,{\em Let $\varphi_p(x)$ be a family of functions that are continuous with respect to $x \!\in\! (a,b)$
for each $p \!\in\! \RR^{+}$ and decreasingly stratified for $p \!\in\! \RR^{+}$, and let $c,d$ be in $\RR^{+}$
such that $c<d$. If$\,:$
\begin{itemize}
\item[{\rm (a)}] $\varphi_c(x) \!>\! 0$ and $\varphi_d(x) \!<\! 0$ for all $x \!\in\! (a,b)$,
    and at the endpoints $\varphi_c(a+) = \varphi_d(b-) = \varphi_d(a+) = 0$ and $\varphi_c(b-) \!\in\! \RR^{+}$ hold;
\item[{\rm (b)}] the functions $\varphi_p(x)$ are continuous
    with respect to $p \!\in\! (c,d)$ for each $x \!\in\! (a,b)$
    and $\varphi_p(b-)$ is continuous with respect to $p \!\in\! (c,d)$ too;
\item[{\rm (c)}] for all $p \!\in\! (c,d)$, there exists a right neighbourhood of the point $a$ in which
    \mbox{$\varphi_p(x) \!<\! 0$} holds;
\item[{\rm (d)}] for all $p \!\in\! (c,d)$  the function $\varphi_p(x)$ has exactly one extremum $t^{(p)}$ on $(a,b)$, which is minimum;
\end{itemize}
then there exists exactly one solution $p_0$, for $p \!\in\! \RR^{+}$, of the following equation
$$
|\varphi_p(t^{(p)})|
=
\varphi_p(b-)
$$
and for $d_0=|\varphi_{p_0}(t^{(p_0)})|=\varphi_{p_0}(b-)$ we have
$$
d_0
\, =
\inf\limits_{\mbox{\scriptsize $p \!\in\! \RR^{+}$}}
\sup\limits_{\mbox{\scriptsize $x \!\in\! (a,b)$}} |\varphi_p(x)|\,.
$$
}

\vspace{-0.3 cm}
\begin{theorem} {\rm (Nike theorem, Theorem 3 \cite{Malesevic_Mihailovic_2021}, Theorem  2.1. \cite{Malesevic_Lutovac_Banjac_2018_2})}
\label{nike_theorem}
Let $f: (0,c) \longrightarrow \RR$ be $m$ times differentiable function
$($for some $m \!\geq\! 2$, $m \!\in\! \NN)$ satisfying the following conditions$:$
\begin{itemize}
\item[{\rm (a)}]
$f^{(m)}(x) \!>\! 0$ for $x \!\in\! (0,c)$;
\item[{\rm (b)}]
there is a right neighbourhood of zero in which the following inequalities are true:
$$
f<0, \, f'<0, \, \ldots, f^{(m-1)}<0;
$$
\item[{\rm (c)}]
there is a left neighbourhood of $c$ in which the following inequalities are true:
$$
f>0, \, f'>0, \, \ldots, f^{(m-1)}>0.
$$
\end{itemize}
Then the function $f$ has exactly one zero $x_0 \!\in\! (0, c)$, and $f(x) \!<\! 0$ for $x \!\in\! (0, x_0)$
and $f(x) \!>\! 0$ for $x \!\in\! (x_0, c)$. Also, the function $f$ has exactly one local minimum on the interval
$(0, c)$. More precisely, there is exactly one point $t \!\in\! (0,c)$ {\big (}in fact $t \!\in\! (0, x_0)${\big )}
such that $f(t) \!<\! 0$ is the smallest value of the function $f$ on the interval $(0, c)$ and particularly on $(0, x_0)$.
\end{theorem}

In a case when it is not possible to apply the Nike theorem, the following theorem is applied, which gives sufficient conditions that the function on the interval has exactly one zero and exactly one minimum (see section 3).

\begin{theorem} {\rm (The Second Nike theorem, Theorem 4 \cite{Malesevic_Mihailovic_2021})}
Let $f\!:\!(0,c) \!\longrightarrow \!\RR$ be $m$ times differentiable function {\big (}$\,$for some \mbox{$m \!\geq\! 2$},
\mbox{$m \!\in\! \NN$}$\,${\big )} satisfying the following conditions$:$

\smallskip
\noindent\;\;
\mbox{\rm (a)}
\begin{minipage}[t]{112.0 mm}
$f^{(m)}$ has exactly one zero $x_m$ on $(0,c)$ such that $f^{(m)} \!>\! 0$ on $(0,x_m)$ and $f^{(m)}<0$ on
$(x_m,c)$;
\end{minipage}

\smallskip
\noindent\;\;
\mbox{\rm (b)}
\begin{minipage}[t]{112.0 mm}
there is a right neighbourhood of zero in which the following inequalities are true:
\end{minipage}
$$
f<0, f'<0,..., f^{(m-1)}<0;
$$

\smallskip
\noindent\;\;
\mbox{\rm (c)}
\begin{minipage}[t]{112.0 mm}
there is a left neighbourhood of $c$ in which the following inequalities are true:
\end{minipage}
$$
f>0, f'>0,..., f^{(m-1)}>0.
$$
Then the function $f$ has exactly one zero $x_0 \!\in\! (0,c)$ and $f(x) \!<\! 0$ for $x \!\in\! (0,x_0)$
and $f(x) \!>\! 0$ for $x \!\in\! (x_0,c)$. The function $f$ has exactly one minimum on the interval
$(0, c)$, i.e. there is exactly one point $t \!\in\! (0,c)$ {\big (}in fact $t \!\in\! (0, x_0)${\big )}
such that $f(t) \!<\! 0$ is the smallest value of the function $f$ on the interval $(0, c)$ and particularly
on $(0, x_0)$.
\end{theorem}

\begin{remark}
Let us emphasize that the previous two forms of the Nike theorem ensure the existence of a minimax approximant. Also, these two theorems claim that the local minimum at $t$ is the only extremum of the function $f$ on $(0,c)$, which is shown in their proofs $($see {\rm $\cite{Malesevic_Mihailovic_2021}$} \!\!$)$.
\end{remark}

The topic of stratification appears in the recently published papers \cite{Malesevic_Micovic_2023}-\cite{Malesevic_Mihailovic_NenezicJovic_Milinkovic_2022}.

\bigskip\noindent
{\bf 2.2. A METHOD FOR PROVING MTP INEQUALITIES}

\medskip

MTP -- Mixed Trigonometric Polynomial function is determined by:
\begin{equation*}
\label{MTP_function}
f(x) = \sum^n_{i=1} \alpha_i x^{p_i} \cos^{q_i}\!x \sin^{r_i}\!x \,,
\end{equation*}
where $x \!\in\! \SSS \!\subseteq\! \RR$ ($\SSS$ is an open or closed interval),
$\alpha_i \!\in\! \RR \!\setminus\! \{0\}, p_i,q_i,r_i \!\in\! \NN_0$ and $n \!\in\! \NN$.
The corresponding inequality
\begin{equation*}
\label{MTP_inequality}
f(x) > 0 \,,
\end{equation*}
where $x \in \SSS$,
is called an MTP inequality.

MTP functions were originally considered through MTP systems of such functions involving multiple variables, see \cite{Dong_Yu_2008}. In the article \cite{Malesevic_Makragic_2016}, the previous definitions of MTP function and MTP inequality were introduced. Moreover, in that article, a method for proving MTP inequalities on the base interval $(0, \pi/2)$ was presented. The method is based on Maclaurin approximations of the sine and cosine functions.
Let us emphasize that in \cite{Shiping_Zhong_2016}, a method for proving MTP inequalities on the base interval, based on the universal trigonometric substitution and Maclaurin approximations, was presented, see also \cite{Chen_Liu_2020}. It is noteworthy that both methods have been applied earlier in numerous papers and monographs, see, for example, \cite{Mortici_2011}, \cite{Mitrinovic_1970}, \cite{Milovanovic_Rassias_2014}. The topic of MTP functions and stratified families of functions has been the subject of recently defended doctoral dissertations \cite{Makragic_2018}, \cite{Banjac_2019}, \cite{Nenezic_2023}.

In further consideration, let $f(x)$ be an MTP function with rational coefficients.
The method for proving MTP inequalities from the paper \cite{Malesevic_Makragic_2016} has been computationally implemented through the doctoral dissertation \cite{Banjac_2019}. For methods for proving inequalities by computer, see also the papers \cite{Malesevic_2007} and \cite{Banjac_Makragic_Malesevic_2016}.
The computer implementation from \cite{Banjac_2019}, called SimTheP, proves MTP inequalities on the interval $\SSS \!\subseteq\! [0, \pi/2]$, providing users with the proof in four stages. Therefore, we will outline the method through a brief description of each of these stages based on the papers
\cite{Malesevic_Makragic_2016}, \cite{Malesevic_Jovanovic_2024}, \cite{Malesevic_Banjac_2019}.

\break
\noindent
\underline{\textbf{I Recognition of possible case}}

\smallskip

In the first phase, the hypothesis $f(x) > 0$ on $\SSS$ is tested
based on the values of the function at the boundary points of the interval.

\medskip
\noindent
\underline{\textbf{II Transformation of angles}}

\smallskip

In this phase, each addend of the MTP function $f(x)$ undergoes the substitution of the expression $\cos^n \! x \sin^m \! x$ $(n, m \!\in\! \NN_0)$
into a sum of sine and cosine functions of multiple angles according to Table 1. These substitutions were proven in the paper \cite{Malesevic_Makragic_2016}.

\vspace*{-2.5 mm}

\begin{table}[H]
\caption{Substitutions of terms $x^{p_i} \cos^{q_i}\!x \sin^{r_i}\!x$, see Table 1 from \cite{Malesevic_Banjac_2019}}

\vspace*{0.25 mm}

\centering{}%
  \begin{tabular}{|c|c|c|}
    \hline
    \multicolumn{3}{|c|}{$\mbox{\boldmath $\cos$}^{n}x \, \mbox{\boldmath $\sin$}^{m}x$}\\
    \hline
    $n = q_i$ & $m = r_i$ &  Substitution \tabularnewline
\hline
\hline
    even & even &  $\begin{array}{l} \vspace{-0.35cm} \\
                \mathop{\mbox{\large $\sum$}}\limits_{k=0}^{\frac{n}{2}+\frac{m}{2}-1}
                \mathop{\mbox{\large $\sum$}}\limits_{j=0}^{k}
                    \dfrac{ (-1)^{\frac{m}{2}+k+j}{n \choose j}{m \choose k-j}
                    \mbox{\boldmath  $\cos$} \left((n+m-2k)x\right)}{2^{n+m-1}}\\[3.0 ex]
                +\;\mathop{\mbox{\large $\sum$}}\limits_{j=0}^{\frac{n}{2}+\frac{m}{2}}
                    \dfrac{ (-1)^{m+\frac{n}{2}+j}{n \choose j}{m \choose \frac{n}{2}+\frac{m}{2}-j}}{2^{n+m}}
                \end{array}$\\
    \hline
    odd & even & $\begin{array}{l} \vspace{-0.35cm} \\
                \mathop{\mbox{\large $\sum$}}\limits_{k=0}^{\frac{n}{2}+\frac{m}{2}-\frac{1}{2}}
                \mathop{\mbox{\large $\sum$}}\limits_{j=0}^{k}
                    \dfrac{(-1)^{\frac{m}{2}+k+j} {n\choose j} {m \choose {k-j}}
                    \mbox{\boldmath $\cos$} \left(\left(n+m-2k\right)x\right)}{2^{n+m-1}}\end{array}$\\
    \hline
    even & odd & $\begin{array}{l} \vspace{-0.35cm} \\
        \mathop{\mbox{\large $\sum$}}\limits_{k=0}^{\frac{n}{2}+\frac{m}{2}-\frac{1}{2}}
        \mathop{\mbox{\large $\sum$}}\limits_{j=0}^{k}
            \dfrac{(-1)^{\frac{m}{2}+k+j-\frac{1}{2}} {n\choose j} {m \choose {k-j}}
            \mbox{\boldmath $\sin$} \left(\left(n+m-2k\right)x\right)}{2^{n+m-1}}

    \end{array}$\\
    \hline
    odd & odd & $\begin{array}{l} \vspace{-0.35cm} \\
        \mathop{\mbox{\large $\sum$}}\limits_{k=0}^{\frac{n}{2}+\frac{m}{2}-1}
        \mathop{\mbox{\large $\sum$}}\limits_{j=0}^{k}
            \dfrac{(-1)^{\frac{m}{2}+k+j-\frac{1}{2}} {n\choose j} {m \choose {k-j}}
            \mbox{\boldmath $\sin$} \left(\left(n+m-2k\right)x\right)}{2^{n+m-1}}
    \end{array}$\\
    \hline
  \end{tabular}
\end{table}

The initial MTP function is transformed into the equivalent form:
\begin{equation}
\label{MTP_function_2}
f(x)
=
\mathop{\mbox{$\displaystyle\sum$}}\limits_{i=1}^{n} \alpha_{i} \, x^{p_i} \! \left(\,
\mathop{\mbox{$\displaystyle\sum$}}\limits_{k=0}^{m_i}
\theta_k \, \mbox{\bf trig}_{k}^{(q_i,r_i)} \! {\Big (}
\underbrace{\mathop{\left(\mbox{\small $q_i - r_i - 2k$}\right) x}}\limits_{(=t)} {\Big )} \!\right) ,
\end{equation}
\noindent
where
$$
\mbox{\bf trig}_{k}^{(q_i,r_i)}
=
\left\{
\begin{array}{ccc}
\mbox{\boldmath $\cos$} \!&\!:\!&\! \mbox{\footnotesize $q_i\mbox{-odd}, r_i\mbox{-even} \;\mbox{or} \; q_i\mbox{-even}, r_i\mbox{-even}$} \\[0.35 ex]
\mbox{\boldmath $\sin$} \!&\!:\!&\! \mbox{\footnotesize $q_i\mbox{-odd}, r_i\mbox{-odd} \;\mbox{or} \; q_i\mbox{-even}, r_i\mbox{-odd}$}
\end{array}
\right.
$$
and
$$
m_i = m_i(\mbox{\small $q_i$}, \mbox{\small $r_i$}) = \left\lceil \frac{q_i + r_i}{2} \right\rceil - 1.
$$

\break

\medskip
\noindent
\underline{\textbf{III Determination of downward rational polynomial approximation}}

\smallskip

During this phase, it is necessary to determine a downward polynomial approximation with rational coefficients of the MTP function (\ref{MTP_function_2}).

Let us specify the general concept of downward/upward polynomial approximation of a function.
Let $\phi(x) : \AAA \longrightarrow \RR$ be any function defined over $\AAA \subseteq \RR$.
The downward polynomial approximation of the function $\phi(x)$ over $\BBB \subseteq \AAA$
is a polynomial $P(x)$ such that
$$
\left(\forall x \in \BBB \right) \phi(x) \geq P(x) \,,
$$
denoted by $\underline{P}(x)$.
Similarly, the upward polynomial approximation of the function $\phi(x)$ over $\BBB \subseteq \AAA$
is a polynomial $P(x)$ such that
$$
\left(\forall x \in \BBB \right) \phi(x) \leq P(x) \,,
$$
denoted by $\overline{P}(x)$.

In the following Lemma from the paper \cite{Malesevic_Jovanovic_2024}, we provide some upward and downward polynomial approximations of the sine and cosine functions.  These assertions were proven in the paper \cite{Malesevic_Makragic_2016}.

\begin{lemma}
\label{mtp_tejlor} {\rm (Lemma 1 \cite{Malesevic_Jovanovic_2024})}
It holds$:$

\medskip\noindent
{\boldmath $(a)$} For the polynomial
$$
T_n(t) = \sum^{(n-1)/2}_{i=0} \dfrac{(-1)^i t^{2i+1}}{(2i+1)!} \,,
$$
where $n = 4k + 1, k \!\in\! \NN_0 \,,$ it holds$:$
\begin{equation*}
\left(\forall t \in \left[ 0, \sqrt{(n+3)(n+4)} \right] \right) \overline{T}_n(t) \geq \overline{T}_{n+4}(t) \geq \sin t \,,
\end{equation*}
\begin{equation*}
\left(\forall t \in \left[ - \sqrt{(n+3)(n+4)}, 0 \right] \right) \underline{T}_n(t) \leq \underline{T}_{n+4}(t) \leq \sin t \,.
\end{equation*}
For $t=0$, the inequalities turn into equalities.
For $t = \pm \sqrt{(n+3)(n+4)}$, the equalities $\overline{T}_n(t) = \overline{T}_{n+4}(t)$
and $\underline{T}_n(t) = \underline{T}_{n+4}(t)$ hold, respectively.

\medskip\noindent
{\boldmath $(b)$} For the polynomial
$$
T_n(t) = \sum^{(n-1)/2}_{i=0} \dfrac{(-1)^i t^{2i+1}}{(2i+1)!} \,,
$$
where $n = 4k + 3, k \!\in\! \NN_0 \,,$ it holds$:$
\begin{equation*}
\left(\forall t \in \left[ 0, \sqrt{(n+3)(n+4)} \right] \right) \underline{T}_n(t) \leq \underline{T}_{n+4}(t) \leq \sin t \,,
\end{equation*}
\begin{equation*}
\left(\forall t \in \left[ - \sqrt{(n+3)(n+4)}, 0 \right] \right) \overline{T}_n(t) \geq \overline{T}_{n+4}(t) \geq \sin t \,.
\end{equation*}
For $t=0$, the inequalities turn into equalities.
For $t = \pm \sqrt{(n+3)(n+4)}$, the equalities $\underline{T}_n(t) = \underline{T}_{n+4}(t)$
and $\overline{T}_n(t) = \overline{T}_{n+4}(t)$ hold, respectively.

\medskip\noindent
{\boldmath $(c)$} For the polynomial
$$
T_n(t) = \sum^{n/2}_{i=0} \dfrac{(-1)^i t^{2i}}{(2i)!} \,,
$$
where $n = 4k, k \!\in\! \NN_0 \,,$ it holds$:$
\begin{equation*}
\left(\forall t \in \left[ -\sqrt{(n+3)(n+4)}, \sqrt{(n+3)(n+4)} \right] \right) \overline{T}_n(t)
\geq \overline{T}_{n+4}(t) \geq \cos t \,.
\end{equation*}
For $t=0$, the inequalities turn into equalities.
For $t = \pm \sqrt{(n+3)(n+4)}$, the equality $\overline{T}_n(t) = \overline{T}_{n+4}(t)$ holds.

\medskip\noindent
{\boldmath $(d)$} For the polynomial
$$
T_n(t) = \sum^{n/2}_{i=0} \dfrac{(-1)^i t^{2i}}{(2i)!} \,,
$$
where $n = 4k+2, k \!\in\! \NN_0 \,,$ it holds$:$
\begin{equation*}
\left(\forall t \in \left[ -\sqrt{(n+3)(n+4)}, \sqrt{(n+3)(n+4)} \right] \right) \underline{T}_n(t)
\leq \underline{T}_{n+4}(t) \leq \cos t \,.
\end{equation*}
For $t=0$, the inequalities turn into equalities.
For $t = \pm \sqrt{(n+3)(n+4)}$, the equality $\underline{T}_n(t) = \underline{T}_{n+4}(t)$ holds.
\end{lemma}

Let ${T}_n^{\phi,a}(x)$ denote the Taylor expansion of order $n$ of some analytic function $\phi$ in the neighbourhood of some point $a$.

With the aim of obtaining a downward polynomial approximation $P(x)$ of the MTP function $f(x)$, we approximate each addend of the function (\ref{MTP_function_2}) by a Maclaurin polynomial using the following estimates:
$$
\left\{
\begin{array}{cc}
\alpha_i \theta_k > 0 : & \mbox{\boldmath $\cos$} \, t > T^{\cos, 0}_{4\ell_1+2}(x),                   \\[1.0 ex]
\alpha_i \theta_k < 0 : & \mbox{\boldmath $\cos$} \, t < T^{\cos, 0}_{4\ell_2}(x),                     \\[1.0 ex]
\alpha_i \theta_k > 0 : & \mbox{\boldmath $\sin$} \, t > T^{\sin, 0}_{4\ell_3+3}(x),                   \\[1.0 ex]
\alpha_i \theta_k < 0 : & \mbox{\boldmath $\sin$} \, t < T^{\sin, 0}_{4\ell_4+1}(x);
\end{array}
\right.
\leqno (\ast)
$$
where $t = \left(q_i-r_i-2k\right) x$ and $\ell_{1, 2, 3, 4} \!\in\! \NN_0$.

By applying $(\ast)$, we determine a polynomial $P(x)$ such that
$$
f(x) > P(x)
$$
for $x \!\in\!  \SSS$.
If there exists a polynomial $P(x)$ with rational coefficients such that
$$
P(x) > 0
$$
for $x \!\in\!  \SSS$, then
$$
f(x) > 0
$$
for $x \!\in\!  \SSS$.

\break

If the coefficients of the MTP function $f(x)$ are not rational numbers but computable real numbers, then we could determine a downward polynomial approximation with rational numbers, see the paper \cite{Malesevic_Banjac_2020}.

\medskip
\noindent
\underline{\textbf{IV The final part}}

\smallskip

For real polynomials defined over a segment with endpoints where the polynomial does not have zero, Sturm's theorem provides the number of roots on such a segment,
see, for example, Theorem {\rm 4.1} \cite{Cutland_1980} or originally
\cite{Sturm_1829}. It is particularly noteworthy that for {\em polynomial functions with rational coefficients defined over a segment with rational endpoints}, according to
Theorem {\rm 4.2} \cite{Cutland_1980}, the problem of determining the number of roots over that segment is
{\em an algorithmically decidable problem}. For such polynomial functions, if we obtain a proof of positivity using Sturm's theorem, we can consider it as
{\em an effective proof} by finite procedures that can be manually verified.

In the third part, $P(x)$ is determined as a polynomial with rational coefficients. If $\SSS$ is not a segment with rational endpoints or the polynomial $P(x)$ has a root at the boundary points of the segment $\SSS$, we consider the polynomial over an extended segment with rational endpoints, see \cite{Malesevic_Banjac_SesumCavic_Korolija_2022}. It is always possible to choose such a segment with rational endpoints that the polynomial $P(x)$  does not have a root at the boundary points of that segment. If the number of roots does not increase over such an extended segment, and we know whether the polynomial has a root at the boundary points of the segment $\SSS$,  then we also have an effective proof of the polynomial inequality $P(x) > 0$ over $\SSS$ by applying Sturm's theorem.

Let us emphasize that proofs based on the Sturm algorithm are absolutely theoretically rigorous, as pointed out in the paper \cite{Alzer_Kwong_2016}, see also
\cite{Alzer_Kwong_2017}-\cite{Alzer_Kwong_2024}.

In  Appendix sections {\bf A1} and {\bf A4} - {\bf A7}, we prove MTP inequalities $f(x) > 0$ over the base interval $(0, \pi/2)$, while in section {\bf A2}, we prove MTP inequality over the interval $(0, 1]$ and in {\bf A3} over the interval $[1, \pi/2]$.

\vspace{1.5cc}
\begin{section}
{MAIN RESULTS}
\end{section}

According to the paper \cite{Malesevic_Mihailovic_2021}, the following statements, which are improvements of Theorems 1--6 from the paper \cite{Mortici_2011}, are proved. Note that the automatic prover SimTheP was utilized for proving the MTP inequalities. The results obtained by this prover are provided in the Appendix.

% Statement 1 Statement 1 Statement 1 Statement 1 Statement 1 Statement 1 Statement 1 Statement 1 Statement 1 Statement 1
% ***********************************************************************************************************************

\bigskip\noindent
\underline{\textbf{Improvement of Theorem 1}}
\begin{lemma}
The family of functions
$$
\hspace*{25 mm}
\varphi_p(x)
=
-\cos x \, + \, \left(\dfrac{\sin x}{x}\right)^{\!3} - \, \dfrac{1}{15} \, x^4 \, + \, p \, x^6
\hspace*{8 mm} {\Big (}\mbox{for $x \in \left(0, \dfrac{\pi}{2} \right)$}{\Big )}
$$
is increasingly stratified with respect to parameter $p \!\in\! \RR^{+}$.
\end{lemma}

Let us introduce the function $g(x)$ so that the equivalence
$$
\varphi_{p} (x)=0 \Longleftrightarrow p=g(x), \,\,\,\,\, x\in \left(0, \dfrac{\pi}{2} \right), \,p\in \mathbb{R^{+}}
$$
holds. Then
$$
g(x) = \dfrac{x^7 + 15\, x^3 \cos x  - 15 \sin^3 x}{15 x^9},
\hspace*{5 mm} x \in \left(0, \dfrac{\pi}{2} \right) .
$$
Note that
$$
\varphi_{p} (x)=\big{(}p-g(x)\big{)}x^{6} \,.
$$

\begin{lemma}
The function $g(x)$ is strictly decreasing for $x\in (0, \pi/2)$.
\end{lemma}
\begin{proof}
Let us notice that the derivative $g'$ is
$$
\!\!
\begin{array}{rcl}
g'(x)
\!\!\!&\!\!=\!\!&\!\!\!
\dfrac{
45 x \cos^3\!x
\!-\!
135 \cos^2\!x \sin x
\!-\!
\left(90 x^3 \!+\! 45 x\right) \cos x
\!-\!
15  x^4 \sin x \!+\! 135 \sin x \!-\! 2 x^7}{15 x^{10}} .
\end{array}
$$
It holds
$$
g'(x)<0, \,\,\,x\in \left(0, \dfrac{\pi}{2} \right) \Longleftrightarrow f(x)>0, \,\,\,x\in \left(0, \dfrac{\pi}{2} \right) ,
$$
where
$$f(x)
=
2 \, x^7
\!-
45 \, x \cos^3\!x
+
135 \cos^2\!x \sin x
+
\left(90 \, x^3 \!+\! 45 \,x\right) \cos x
+
15 \, x^4 \sin x - 135 \sin x.
$$
According to \cite{Malesevic_Makragic_2016}, there exists proof that the MTP function $f$ is positive for $x\!\in \!(0,\pi/2)$. The proof is given in {\bf Appendix A1}.
\end{proof}

\begin{statement}
Let $$
A = \dfrac{4 \pi^7-7680}{15 \pi^9} = 0.0098430 \ldots
\quad\mbox{and}\quad
B = \dfrac{23}{1890} = 0.012169 \ldots \,.
$$
Then, it holds$:$

\smallskip
\noindent
{\boldmath $(i)$}
If $p \in (0,A]$, then
$$
x \in \left(0, \dfrac{\pi}{2} \right)
\Longrightarrow
-\dfrac{1}{15}x^4 + p x^6
\leq
-\dfrac{1}{15}x^4 + A x^6
<
\cos x - \left(\dfrac{\sin x}{x}\right)^{\!3}.
$$

\smallskip
\noindent
{\boldmath $(ii)$}
If $p \in (A,B)$, then $\varphi_{p}(x)$ has exactly one zero $x_{0}^{(p)}$ on $(0, \pi/2)$. Also,
$$
x \in \left(0,x_{0}^{(p)} \right)
\Longrightarrow
-\dfrac{1}{15}x^4 + p \, x^6
<
\cos x - \left(\dfrac{\sin x}{x}\right)^{\!3}
$$
and

\break

$$
x \in \left(x_{0}^{(p)},\dfrac{\pi}{2} \right)
\Longrightarrow
\cos x - \left(\dfrac{\sin x}{x}\right)^{\!3}
<
-\dfrac{1}{15}x^4 + p \, x^6
$$
hold.

\smallskip
\noindent
{\boldmath $(iii)$}
If $p \in [B,\infty)$, then
$$
x \in \left(0, \dfrac{\pi}{2} \right)
\Longrightarrow
\cos x - \left(\dfrac{\sin x}{x}\right)^{\!3}
<
-\dfrac{1}{15}x^4 + B x^6
\leq
-\dfrac{1}{15}x^4 + p x^6.
$$

\smallskip
\noindent
{\boldmath $(iv)$}
There is exactly one solution to the following equation
$$
\left| \varphi_p{\big (}t^{(p)}{\big )}\right| =  \varphi_p \left(\frac{\pi}{2}- \right) ,
$$
where $t^{(p)}$ is a minimum of  $\varphi_p(x)$ on $(0, \pi/2)$, with respect to parameter $p \!\in\! (A,B)$, which is numerically determined as
$$
p_0 = 0.010004 \ldots \,.
$$
For the value
$$
d_{0}
=
\left| \varphi_{p_0}{\big (}t^{(p_0)}{\big )} \right|=  \varphi_{p_0} \left(\frac{\pi}{2}- \right)
=
0.0024209\ldots \,,
$$
the following holds
$$
d_0 = \inf\limits_{p \in \RR^{+}} \sup\limits_{x \in (0,\pi/2)}{\left| \varphi_p(x) \right|}.
$$

\smallskip
\noindent
{\boldmath $(v)$}
For the value $p_0 = 0.010004 \ldots  $ , the minimax approximant of the family is determined as
$$
\varphi_{p_0}(x)
=
-\cos x + \left(\dfrac{\sin x}{x}\right)^{\!3}
-
\dfrac{1}{15}x^4
+
0.010004 \ldots   x^6
$$
and it determines the corresponding minimax approximation
\begin{equation*}
\label{minimax_approximation_T1}
\cos x - \left(\dfrac{\sin x}{x}\right)^{\!3}
\approx
-\dfrac{1}{15}x^4 + 0.010004 \ldots   x^6 \,.
\end{equation*}
\end{statement}
\begin{proof}
This statement is based on the results of the paper \cite{Mortici_2011}
and the fact that
$$
A \;= \lim\limits_{x \rightarrow \pi/2-}{g(x)}
\quad\mbox{and}\quad
B \;=\!\! \lim\limits_{x \rightarrow 0+}{g(x)} \,.
$$

\noindent
The function $g(x)$ is continuous and, according to Lemma 3, strictly decreasing on $(0, \pi/2)$.

\break

\smallskip
\noindent
{\boldmath $(i)$}
If $p\in (0,A]$, then
$$
g(x)>A \Longleftrightarrow (A-g(x))  x^{6}<0 \Longleftrightarrow \varphi_{A}(x)<0
$$
and, therefore, we can conclude that
$$
x \in \left(0, \dfrac{\pi}{2} \right)
\Longrightarrow
-\dfrac{1}{15}x^4 + p x^6
\leq
-\dfrac{1}{15}x^4 + A x^6
<
\cos x - \left(\dfrac{\sin x}{x}\right)^{\!3}.
$$

\smallskip
\noindent
{\boldmath $(ii)$}
If $p\in (A,B)$, based on Lemma 3, the equation
$$
g(x)=p
$$
has a unique solution $x_{0}^{(p)}$ and it holds
$$
\begin{array}{rcl}
x \in \left(0,x_{0}^{(p)} \right)
\Longrightarrow
g(x)>p
\!\!&\!\!\Longleftrightarrow\!\!&\!\!
\varphi_{p}(x)<0                             \\[3.0 ex]
\!\!&\!\! \Longleftrightarrow\!\!&\!\!
-\dfrac{1}{15}x^4 + p \, x^6
<
\cos x - \left(\dfrac{\sin x}{x}\right)^{\!3}
\end{array}
$$
and
$$
\begin{array}{rcl}
x \in \left(x_{0}^{(p)}, \dfrac{\pi}{2} \right)
\Longrightarrow
g(x)<p
\!\!&\!\!\Longleftrightarrow\!\!&\!\!
\varphi_{p}(x)>0                             \\[3.0 ex]
\!\!&\!\! \Longleftrightarrow\!\!&\!\!
\cos x - \left(\dfrac{\sin x}{x}\right)^{\!3}
<
-\dfrac{1}{15}x^4 + p \, x^6.
\end{array}
$$

\smallskip
\noindent
{\boldmath $(iii)$}
If $p\in [B,\infty)$, then
$$
g(x)<B \Longleftrightarrow (B-g(x)) x^{6}>0 \Longleftrightarrow \varphi_{B}(x)>0
$$
and, therefore, we can conclude that
$$
x \in \left(0, \dfrac{\pi}{2} \right)
\Longrightarrow
\cos x - \left(\dfrac{\sin x}{x}\right)^{\!3}
<
-\dfrac{1}{15}x^4 + B x^6
\leq
-\dfrac{1}{15}x^4 + p x^6.
$$

\smallskip
\noindent
{\boldmath $(iv), (v)$}
Let $p\in (A,B)$. For the family $\varphi_p(x)$, the Taylor's expansions are:
\begin{equation}
\label{Tejlor_11}
\varphi_p(x)
=
\left(- \dfrac{23}{1890} + p\right) x^6 + \dfrac{41}{37800} x^8 + o(x^8)
\end{equation}
and
\begin{equation}
\begin{array}{rcl}
\label{Tejlor_21}
\varphi_p(x)
\!\!&\!\!=\!\!&\!\!
\left( \dfrac{8}{\pi^3} - \dfrac{\pi^4}{240} + \dfrac{\pi^6}{64} \,p \right) +                    \\[3.0 ex]
\!\!&\!\! \!\!&\!\!
\!\!\!\!\!\!+
\left( -\dfrac{48}{\pi^4} + 1 - \dfrac{\pi^3}{30} + \dfrac{3\, \pi^5}{16} \, p \right)
\left(x-\dfrac{\pi}{2}\right) +                                                                   \\[3.0 ex]
\!\!&\!\! \!\!&\!\!
\!\!\!\!\!\!+
\left( \dfrac{192}{\pi^5} - \dfrac{12}{\pi^3} - \dfrac{\pi^2}{10} + \dfrac{15 \, \pi^4}{16} \, p \right)
\left(x-\dfrac{\pi}{2}\right)^{\!2} +                                                             \\[3.0 ex]
\!\!&\!\! \!\!&\!\!
\!\!\!\!\!\!+
\left( -\dfrac{640}{\pi^6} + \dfrac{72}{\pi^4} - \dfrac{1}{6} - \dfrac{2 \, \pi}{15} + \dfrac{5 \, \pi^3}{2} \, p \right)
\left(x-\dfrac{\pi}{2}\right)^{\!3} +                                                              \\[3.0 ex]
\!\!&\!\! \!\!&\!\!
\!\!\!\!\!\!+
\left( \dfrac{1920}{\pi^7} - \dfrac{288}{\pi^5} + \dfrac{7}{\pi^3} - \dfrac{1}{15} + \dfrac{15 \, \pi^2}{4} p \right)
\left(x-\dfrac{\pi}{2}\right)^{\!4} +                                                              \\[3.0 ex]
\!\!&\!\! \!\!&\!\!
\!\!\!\!\!\!+
\left( -\dfrac{5376}{\pi^8} + \dfrac{960}{\pi^6} - \dfrac{42}{\pi^4} + \dfrac{1}{120} + 3 \, \pi \, p\right)
\left(x-\dfrac{\pi}{2}\right)^{\!5} +                                                             \\[3.0 ex]
\!\!&\!\! \!\!&\!\!
\!\!\!\!\!\!+
\left( \dfrac{14336}{\pi^9} - \dfrac{2880}{\pi^7} + \dfrac{168}{\pi^5} - \dfrac{61}{30 \, \pi^3} + p\right)
\left(x-\dfrac{\pi}{2}\right)^{\!6} +                                                             \\[3.0 ex]
\!\!&\!\! \!\!&\!\!
\!\!\!\!\!\!+
\left(
-\dfrac{36864}{\pi^{10}} + \dfrac{8064}{\pi^8} - \dfrac{560}{\pi^6} + \dfrac{61}{5 \, \pi^4} - \dfrac{1}{5040}
\right)
\left(x-\dfrac{\pi}{2}\right)^{\!7} +                                                              \\[3.0 ex]
\!\!&\!\! \!\!&\!\!
\!\!\!\!\!\!+\;
o \left( \left( x - \dfrac{\pi}{2} \right)^{\!7}  \right).
\end{array}
\end{equation}

\noindent
For $p \!\in\! (A,B)$, functions $\varphi_p(x)$ don't satisfy all of the conditions of the Nike theorem. In consequence, we use the Second Nike theorem. Now we check the fulfillment of the Second Nike theorem:

\medskip
\noindent
(a) Let us observe the seventh derivative
$$
\begin{array}{rcl}
\varphi_p^{(7)}(x)
\!\!&\!\!=\!\!&\!\!
\dfrac{1}{x^{10}}
{\Big (}
\left(2187 \, x^7 - 61236 \, x^5 + 340200 \, x^3 - 423360 \, x\right)\cos^3 x \, +                \\[3.0 ex]
\!\!&\!\! \!\!&\!\!
\quad\; +
\left(-15309 \, x^6 + 170100 \, x^4 - 476280 \, x^2 + 181440\right) \sin x \cos^2 x \, +          \\[3.0 ex]
\!\!&\!\! \!\!&\!\!
\quad\; +
\left(-1641 \, x^7 + 46116 \, x^5 - 264600 \, x^3 + 423360 \, x\right) \cos x \, +                \\[3.0 ex]
\!\!&\!\! \!\!&\!\!
\quad\; +
\left(-x^{10} + 3843 \, x^6 - 44100 \, x^4 + 158760 \, x^2 - 181440\right) \sin x
{\Big )},
\end{array}
$$
for $x \!\in\! (0, \pi/2)$. Now we prove that function $\varphi_p^{(7)}(x)$ has exactly one zero $c$ on $(0, \pi/2)$ such that $\varphi_p^{(7)}(x) > 0$ for
$x \!\in\! (0,c)$, and $\varphi_p^{(7)}(x) < 0$ for $x \!\in\! (c, \pi/2)$. Function
$\varphi_p^{(7)}(x)$ is positive for $x \!\in\! [0, 1]$ because, according to \cite{Malesevic_Makragic_2016}, there exists proof that the numerator
of $\varphi_p^{(7)}(x)$ is positive on $[0, 1]$.
The proof is given in {\bf Appendix A2}. Furthermore, let us observe the eighth derivative
$$
\!\!
\begin{array}{rcl}
\varphi_p^{(8)}(x)
\!\!&\!\!=\!\!&\!\!
\dfrac{1}{x^{11}}
{\Big (}
\left(-52488 \, x^7+816480 \, x^5-3810240 \, x^3+4354560 \, x\right) \cos^3 x \, +                \\[3.0 ex]
\!\!&\!\! \!\!&\!\!
\quad\; +
\left(-6561x^8\!\!+\!244944x^6\!\!-\!2041200x^4\!\!+\!5080320x^2\!\!-\!1814400\right)\sin x\cos^2 x\, +\\[3.0 ex]
\!\!&\!\! \!\!&\!\!
\quad\; +
\left(-x^{11} + 39384 \, x^7 - 614880 \, x^5+2963520 \, x^3 - 4354560 \, x\right) \cos x \, +     \\[3.0 ex]
\!\!&\!\! \!\!&\!\!
\quad\; +
\left(1641 \, x^8-61488 \, x^6+529200 \, x^4-1693440 \, x^2+1814400\right) \sin x
{\Big )},
\end{array}
$$
for $x \!\in\! (0, \pi/2)$. According to \cite{Malesevic_Makragic_2016}, there exists proof that the numerator of $\varphi_p^{(8)}(x)$ is negative on $[1, \pi/2]$.
It is enough to prove that $\varphi_p^{(8)}(x)$ is negative on $[1, \pi/2]$ using MacLaurin polynomials.
The proof is given in {\bf Appendix A3}.
Finally,
$$
\varphi_p^{(7)}\!\!\left(\dfrac{\pi}{2}\right)
\!\!=\!\!
\dfrac{
-\pi^{10} \!+\! 61488  \pi^6 \!-\! 2822400  \pi^4 \!+\! 40642560  \pi^2 \!-\! 185794560}{\pi^{10}}
\!\!=\!\!
-6.14 789 \ldots
\!<\!
0\,.
$$
Therefore, there exists exactly one zero $c \!\in\! (0, \pi/2)$ of function
$\varphi_p^{(7)}(x)$ such that $\varphi_p^{(7)}(x) > 0$ for $x \!\in\! (0,c)$ and
$\varphi_p^{(7)}(x) < 0$ for $x \!\in\! (c, \pi/2)$, where $c$ is numerically determined as $c = 1.40749 \ldots$.
It is hereby shown that for $m = 7$, the first condition of the Second Nike theorem is satisfied.

\medskip
\noindent
(b) According to (\ref{Tejlor_11}), there is a right neighbourhood ${\cal U}_{0}$ of zero in which
$$
\hspace*{30.0 mm}
\varphi_p(x)
<
0
\mbox{ ,}\quad
\varphi_p'(x)
<
0
\mbox{ ,}\quad
\ldots
\mbox{ ,}\quad
\varphi_p^{(6)}(x)
<
0
\hspace*{12.0 mm}
{\big (} x \!\in\! {\cal U}_{0} {\big )}
$$ hold.

\medskip
\noindent
(c) According to (\ref{Tejlor_21}), there is a left neighbourhood ${\cal U}_{\pi/2}$ of $\pi/2$ in which
$$
\hspace*{32.0 mm}
\varphi_p(x)
>
0
\mbox{ ,}\quad
\varphi_p'(x)
>
0
\mbox{ ,}\quad
\ldots
\mbox{ ,}\quad
\varphi_p^{(6)}(x)
>
0
\hspace*{8.0 mm}
{\big (} x \!\in\! {\cal U}_{\pi/2} {\big )}
$$ hold.
Then, for every function $\varphi_p(x)$, on the interval $(0, \pi/2)$, there exists exactly one extremum $t$, which is minimum, and there exists exactly one zero $x_0$.
Hence, for the family of functions $\varphi_p(x)$, conditions of the Second Nike theorem are satisfied, as well as conditions of Theorem 7, which implies the existence of a minimax approximant.
Minimax approximant and error can be numerically determined via {\tt Maple} software in the following way: let $f(x,p) 
\, {\tt :=} \, \varphi_p(x)$ and $F(x,p) \, {\tt :=} \, \varphi_p'(x)$, then using the command
{\sf$$
\mbox{\tt fsolve}
{\big (}
\{F(x,p) = 0, \mbox{\tt abs} {\big (} f(x,p) {\big )} = f(\pi/2,p) \}, \{ x = 0 .. \pi/2, p = A .. B \}
{\big )};
$$}
we get
{\sf$$
\{ p = 0.01000418287, x = 1.299862713 \}
$$}
and, for $p_0 = 0.010004 \ldots$ , we get
$$
f(x, p_0)=-\cos x + \left(\dfrac{\sin x}{x}\right)^{\!3}
-
\dfrac{1}{15}x^4
+
0.010004 \ldots x^6
$$
and
$$
d_0
=
f\left( \dfrac{\pi}{2}, p_0 \right) 
=
0.0024209 \ldots \,.
\vspace{-0.3 cm}
$$ 
\end{proof}

Based on the previous considerations, enhancement of Theorem 1 has been obtained in the following form:

\begin{proposition}
For every $0 < x < \pi/2$, the following inequalities hold
$$
-\dfrac{1}{15}x^4 + \dfrac{4 \pi^7-7680}{15 \pi^9} \,x^6
<
\cos x - \left(\dfrac{\sin x}{x}\right)^{\!3}
<
-\dfrac{1}{15}x^4 + \dfrac{23}{1890} \,x^6,
$$
and the constants $A = \dfrac{4 \pi^7-7680}{15 \pi^9} = 0.0098430 \ldots$ and $B = \dfrac{23}{1890} = 0.012169 \ldots$
are the best possible.
\end{proposition}

% Statement 2 Statement 2 Statement 2 Statement 2 Statement 2 Statement 2 Statement 2 Statement 2 Statement 2 Statement 2
% ***********************************************************************************************************************

\smallskip\noindent
\underline{\textbf{Improvement of Theorem 2}}
\begin{lemma}
The family of functions
$$
\hspace*{20.00 mm}
\varphi_p(x)
=
-\dfrac{\sin x}{x} + \dfrac{\cos x + 2}{3} - \dfrac{1}{180}x^4 + p \, x^6
\hspace*{10 mm} {\Big (}\mbox{for $x \in \left(0, \dfrac{\pi}{2} \right)$}{\Big )}
$$
is increasingly stratified with respect to parameter $p \!\in\! \RR^{+}$.
\end{lemma}

Let us introduce the function $g(x)$ so that the equivalence
$$
\varphi_{p} (x)=0 \Longleftrightarrow p=g(x), \,\,\,\,\, x\in \left(0, \dfrac{\pi}{2} \right), \,p\in \mathbb{R^{+}}
$$
holds. Then
$$
g(x) = \dfrac{x^5 - 60 \, x \cos x + 180 \sin x - 120 \, x}{180 \, x^7},
\hspace*{5 mm} x \in \left(0, \dfrac{\pi}{2} \right) .
$$
Note that
$$
\varphi_{p} (x)=\big{(}p-g(x)\big{)} x^{6} \,.
$$

\begin{lemma}
The function $g(x)$ is strictly decreasing for $x\in (0, \pi/2)$.
\end{lemma}
\begin{proof}
Let us notice that the derivative $g'$ is
$$
\begin{array}{rcl}
g'(x)
\!\!\!&\!\!=\!\!&\!\!\!
\dfrac{270 \, x \cos x + 30 \, x^2 \sin x - 630 \sin x - x^5 + 360 \, x }{90 \, x^8}.
\end{array}
$$
It holds
$$
g'(x)<0, \,\,\,x\in \left(0, \dfrac{\pi}{2} \right) \Longleftrightarrow f(x)>0, \,\,\,x\in \left(0, \dfrac{\pi}{2} \right) ,
$$
where
$$f(x)
=
x^5 - 360 \, x - 270 \, x \cos x - 30 \, x^2 \sin x + 630 \sin x.
$$
According to \cite{Malesevic_Makragic_2016}, there exists proof that the MTP function $f$ is positive for $x\!\in \!(0,\pi/2)$. The proof is given in {\bf Appendix A4}.
\end{proof}

\begin{statement}
Let $$
A = \dfrac{\pi^5 \!-\! 1920 \pi \!+\! 5760}{45 \, \pi^7} = 0.00025135 \ldots
\quad\mbox{and}\quad
B = \dfrac{1}{3780} = 0.00026455 \ldots\,.
$$
Then, it holds$:$

\smallskip
\noindent
{\boldmath $(i)$}
If $p \in (0,A]$, then
$$
x \in \left(0, \dfrac{\pi}{2} \right)
\Longrightarrow
- \dfrac{1}{180}x^4 + p \, x^6
\leq
- \dfrac{1}{180}x^4 + A \, x^6
<
\dfrac{\sin x}{x} - \dfrac{\cos x + 2}{3}.
$$

\smallskip
\noindent
{\boldmath $(ii)$}
If $p \in (A,B)$, then $\varphi_{p}(x)$ has exactly one zero $x_{0}^{(p)}$ on $(0, \pi/2)$. Also,
$$
x \in \left(0,x_{0}^{(p)} \right)
\Longrightarrow
- \dfrac{1}{180}x^4 + p \, x^6
<
\dfrac{\sin x}{x} - \dfrac{\cos x + 2}{3}
$$
and
$$
x \in \left(x_{0}^{(p)}, \dfrac{\pi}{2} \right)
\Longrightarrow
\dfrac{\sin x}{x} - \dfrac{\cos x + 2}{3}
<
- \dfrac{1}{180}x^4 + p \, x^6
$$
hold.

\smallskip
\noindent
{\boldmath $(iii)$}
If $p \in [B,\infty)$, then
$$
x \in \left(0, \dfrac{\pi}{2} \right)
\Longrightarrow
\dfrac{\sin x}{x} - \dfrac{\cos x + 2}{3}
<
-\dfrac{1}{180}x^4 + B \, x^6
\leq
-\dfrac{1}{180}x^4 + p \, x^6.
$$

\smallskip
\noindent
{\boldmath $(iv)$}
There is exactly one solution to the following equation
$$
\left| \varphi_p{\big (}t^{(p)}{\big )} \right|=  \varphi_p \left(\frac{\pi}{2}- \right) ,
$$
where $t^{(p)}$ is a minimum of  $\varphi_p(x)$ on $(0, \pi/2)$, with respect to parameter $p \!\in\! (A,B)$, which is numerically determined as
$$
p_0 = 0.00025234 \ldots \,.
$$
For the value
$$
d_{0}
=
\left| \varphi_{p_0}{\big (}t^{(p_0)}{\big )} \right|=  \varphi_{p_0} \left(\frac{\pi}{2}- \right)
=
0.000014887\ldots \,,
$$
the following holds
$$
d_0 = \inf\limits_{p \in \RR^{+}} \sup\limits_{x \in (0,\pi/2)}{\left| \varphi_p(x) \right|}.
$$

\smallskip
\noindent
{\boldmath $(v)$}
For the value $p_0 = 0.00025234 \ldots  $ , the minimax approximant of the family is determined as
$$
\varphi_{p_0}(x)
=
-\dfrac{\sin x}{x} + \dfrac{\cos x + 2}{3} - \dfrac{1}{180}x^4 + 0.00025234 \ldots  x^6
$$
and it determines the corresponding minimax approximation
\begin{equation*}
\label{minimax_approximation_T2}
\dfrac{\sin x}{x} - \dfrac{\cos x + 2}{3} \approx -\dfrac{1}{180}x^4 + 0.00025234 \ldots  x^6 \,.
\end{equation*}
\end{statement}
\begin{proof}
This statement is based on the results of the paper \cite{Mortici_2011} and the fact that
$$
A \;= \lim\limits_{x \rightarrow \pi/2-}{g(x)}
\quad\mbox{and}\quad
B \;=\!\! \lim\limits_{x \rightarrow 0+}{g(x)} \,.
$$

\noindent
The function $g(x)$ is continuous and, according to Lemma 5, strictly decreasing on $(0, \pi/2)$.

\smallskip
\noindent
{\boldmath $(i)$}
If $p\in (0,A]$, then
$$
g(x)>A \Longleftrightarrow (A-g(x)) x^{6}<0 \Longleftrightarrow \varphi_{A}(x)<0
$$
and, therefore, we can conclude that
$$
x \in \left(0, \dfrac{\pi}{2} \right)
\Longrightarrow
- \dfrac{1}{180}x^4 + p \, x^6
\leq
- \dfrac{1}{180}x^4 + A \, x^6
<
\dfrac{\sin x}{x} - \dfrac{\cos x + 2}{3}.
$$

\smallskip
\noindent
{\boldmath $(ii)$}
If $p\in (A,B)$, based on Lemma 5, the equation
$$
g(x)=p
$$
has a unique solution $x_{0}^{(p)}$ and it holds
$$
\begin{array}{rcl}
x \in \left(0,x_{0}^{(p)} \right)
\Longrightarrow
g(x)>p
\!\!&\!\!\Longleftrightarrow\!\!&\!\!
\varphi_{p}(x)<0                             \\[3.0 ex]
\!\!&\!\! \Longleftrightarrow\!\!&\!\!
- \dfrac{1}{180}x^4 + p \, x^6
<
\dfrac{\sin x}{x} - \dfrac{\cos x + 2}{3}
\end{array}
$$
and
$$
\begin{array}{rcl}
x \in \left(x_{0}^{(p)}, \dfrac{\pi}{2} \right)
\Longrightarrow
g(x)<p
\!\!&\!\!\Longleftrightarrow\!\!&\!\!
\varphi_{p}(x)>0                             \\[3.0 ex]
\!\!&\!\! \Longleftrightarrow\!\!&\!\!
\dfrac{\sin x}{x} - \dfrac{\cos x + 2}{3}
<
- \dfrac{1}{180}x^4 + p \, x^6.
\end{array}
$$

\smallskip
\noindent
{\boldmath $(iii)$}
If $p\in [B,\infty)$, then
$$
g(x)<B \Longleftrightarrow (B-g(x)) x^{6}>0 \Longleftrightarrow \varphi_{B}(x)>0
$$
and, therefore, we can conclude that
$$
x \in \left(0, \dfrac{\pi}{2} \right)
\Longrightarrow
\dfrac{\sin x}{x} - \dfrac{\cos x + 2}{3}
<
-\dfrac{1}{180}x^4 + B \, x^6
\leq
-\dfrac{1}{180}x^4 + p \, x^6.
$$

\smallskip
\noindent
{\boldmath $(iv), (v)$}
Let $p\in (A,B)$. For the family $\varphi_p(x)$, the Taylor's expansions are:
\begin{equation}
\label{Tejlor_12}
\varphi_p(x)
=
\left(-\dfrac{1}{3780} + p\right) x^6 + \dfrac{1}{181440} x^8 + o(x^8)
\end{equation}
and
\begin{equation}
\begin{array}{rcl}
\label{Tejlor_22}
\varphi_p(x)
\!\!&\!\!=\!\!&\!\!
\left( -\dfrac{2}{\pi} + \dfrac{2}{3} - \dfrac{\pi^4}{2880}  + \dfrac{\pi^6}{64} \, p \right) +              \\[3.0 ex]
\!\!&\!\! \!\!&\!\!
\!\!\!\!\!\!+
\left( \dfrac{4}{\pi^2} - \dfrac{1}{3} - \dfrac{\pi^3}{360} + \dfrac{3 \, \pi^5}{16} \, p \right)
\left(x-\dfrac{\pi}{2}\right) +                                                                               \\[3.0 ex]
\!\!&\!\! \!\!&\!\!
\!\!\!\!\!\!+
\left( - \dfrac{8}{\pi^3} + \dfrac{1}{\pi} - \dfrac{\pi^2}{120} + \dfrac{15 \, \pi^4}{16} \, p \right)
\left(x-\dfrac{\pi}{2}\right)^{\!2} +                                                                         \\[3.0 ex]
\!\!&\!\! \!\!&\!\!
\!\!\!\!\!\!+
\left( \dfrac{16}{\pi^4} - \dfrac{2}{\pi^2} + \dfrac{1}{18} - \dfrac{\pi}{90} + \dfrac{5 \, \pi^3}{2} \, p \right)
\left(x-\dfrac{\pi}{2}\right)^{\!3} +                                                                         \\[3.0 ex]
\!\!&\!\! \!\!&\!\!
\!\!\!\!\!\!+
\left( - \dfrac{32}{\pi^5} + \dfrac{4}{\pi^3} - \dfrac{1}{12 \, \pi} - \dfrac{1}{180} + \dfrac{15 \, \pi^2}{4} \, p \right)
\left(x-\dfrac{\pi}{2}\right)^{\!4} +                                                                         \\[3.0 ex]
\!\!&\!\! \!\!&\!\!
\!\!\!\!\!\!+
\left( \dfrac{64}{\pi^6} - \dfrac{8}{\pi^4} + \dfrac{1}{6 \, \pi^2} - \dfrac{1}{360} + 3 \, \pi \, p \right)
\left(x-\dfrac{\pi}{2}\right)^{\!5} +                                                                         \\[3.0 ex]
\!\!&\!\! \!\!&\!\!
\!\!\!\!\!\!+
\left( - \dfrac{128}{\pi^7} + \dfrac{16}{\pi^5} - \dfrac{1}{3 \, \pi^3} + \dfrac{1}{360 \, \pi} + p \right)
\left(x-\dfrac{\pi}{2}\right)^{\!6} +                                                                         \\[3.0 ex]
\!\!&\!\! \!\!&\!\!
\!\!\!\!\!\!+
\left( \dfrac{256}{\pi^8} - \dfrac{32}{\pi^6} + \dfrac{2}{3 \, \pi^4} - \dfrac{1}{180 \, \pi^2} + \dfrac{1}{15120} \right)
\left(x-\dfrac{\pi}{2}\right)^{\!7} +                                                                          \\[3.0 ex]
\!\!&\!\! \!\!&\!\!
\!\!\!\!\!\! + \;
o \left( \left( x - \dfrac{\pi}{2} \right)^{\!7}  \right).
\end{array}
\end{equation}
For $p \!\in\! (A,B)$, functions $\varphi_p(x)$ satisfy all of the conditions of the Nike theorem. Now we check the fulfillment of the Nike theorem:

\medskip
\noindent
(a) Let us observe the seventh derivative
$$
\begin{array}{rcl}
\varphi_p^{(7)}(x)
\!\!&\!\!=\!\!&\!\!
\dfrac{1}{3x^8}
{\Big (}
\left(x^8 -21 \, x^6 + 630 \, x^4 - 7560 \, x^2 + 15120 \right) \sin x \, +                       \\[3.0 ex]
\!\!&\!\! \!\!& + \,
3 \, x \left( x^6 - 42 \, x^4 + 840 \, x^2 - 5040\right) \cos x \, {\Big )}
\end{array}
$$
for $x \!\in\! (0, \pi/2)$. According to \cite{Malesevic_Makragic_2016}, there exists proof that the numerator of $\varphi_p^{(7)}(x)$ is positive on $[0, \pi/2]$. The proof is given in {\bf Appendix A5}.
Thus, for $m = 7$, the first condition of the Nike theorem is satisfied.

\medskip
\noindent
(b) According to (\ref{Tejlor_12}), there is a right neighbourhood ${\cal U}_{0}$ of zero in which
$$
\hspace*{30.0 mm}
\varphi_p(x)
<
0
\mbox{ ,}\quad
\varphi_p'(x)
<
0
\mbox{ ,}\quad
\ldots
\mbox{ ,}\quad
\varphi_p^{(6)}(x)
<
0
\hspace*{12.0 mm}
{\big (} x \!\in\! {\cal U}_{0} {\big )}
$$ hold.

\medskip
\noindent
(c) According to (\ref{Tejlor_22}), there is a left neighbourhood ${\cal U}_{\pi/2}$ of $\pi/2$ in which
$$
\hspace*{32.0 mm}
\varphi_p(x)
>
0
\mbox{ ,}\quad
\varphi_p'(x)
>
0
\mbox{ ,}\quad
\ldots
\mbox{ ,}\quad
\varphi_p^{(6)}(x)
>
0
\hspace*{8.0 mm}
{\big (} x \!\in\! {\cal U}_{\pi/2} {\big )}
$$ hold.
Then, for every function $\varphi_p(x)$, on the interval $(0,\pi/2)$, there exists exactly one extremum $t$, which is minimum, and there exists exactly one zero $x_0$.
Hence, for the family of functions $\varphi_p(x)$, conditions of the Nike theorem are satisfied, as well as conditions of Theorem 7, which implies the existence of a minimax approximant.
Minimax approximant and error can be numerically determined via {\tt Maple} software in the following way: 
$f(x,p) \, {\tt :=} \, \varphi_p(x)$ and $F(x,p) \, {\tt :=} \, \varphi_p'(x)$, then using the command
{\sf$$
\mbox{\tt fsolve}
{\big (}
\{F(x,p) = 0, \mbox{\tt abs} {\big (} f(x,p) {\big )} = f(\pi/2,p) \}, \{ x = 0 .. \pi/2, p = A .. B \}
{\big )};
$$}
we get
{\sf$$
\{ p = 0.000252341144, x = 1.305655179 \}
$$}
and, for $p_0 = 0.00025234 \ldots$ , we get
$$
f(x, p_0)=
-\dfrac{\sin x}{x} + \dfrac{\cos x + 2}{3} - \dfrac{1}{180}x^4 + 0.00025234 \ldots x^6
$$
and
$$
d_0
=
f\left( \dfrac{\pi}{2}, p_0 \right) 
=
0.000014887 \ldots \,.
\vspace{-0.3 cm}
$$
\end{proof}

Based on the previous considerations, enhancement of Theorem 2 has been obtained in the following form:

\begin{proposition}
For every $0 < x < \pi/2$, the following inequalities hold
$$
-\dfrac{1}{180}x^4 + \dfrac{\pi^5 - 1920 \pi + 5760}{45 \, \pi^7} \,x^6
<
\dfrac{\sin x}{x} - \dfrac{\cos x + 2}{3}
<
-\dfrac{1}{180}x^4 + \dfrac{1}{3780} \,x^6,
$$
and the constants $A = \dfrac{\pi^5 - 1920 \pi + 5760}{45 \, \pi^7} = 0.00025135 \ldots$
and $B = \dfrac{1}{3780} = 0.00026 455 \ldots \, $ are the best possible.
\end{proposition}

% Statement 3 Statement 3 Statement 3 Statement 3 Statement 3 Statement 3 Statement 3 Statement 3 Statement 3 Statement 3
% ***********************************************************************************************************************

\smallskip\noindent
\underline{\textbf{Improvement of Theorem 3}}
\begin{lemma}
The family of functions
$$
\hspace*{5.00 mm}
\varphi_p(x)
=
2 \, \dfrac{\sin x}{x} + \dfrac{\tan x}{x} - 3 - \dfrac{3}{20} \, x^4 \! \dfrac{1}{\cos x} + p \, x^6 \dfrac{1}{\cos x}
\hspace*{2.25 mm} {\Big (}\mbox{for $x \in \left(0, \dfrac{\pi}{2} \right)$}{\Big )}
$$
is increasingly stratified with respect to parameter $p \!\in\! \RR^{+}$.
\end{lemma}

Let us introduce the function $g(x)$ so that the equivalence
$$
\varphi_{p} (x)=0 \Longleftrightarrow p=g(x), \,\,\,\,\, x\in\left(0, \dfrac{\pi}{2} \right), \,p\in \mathbb{R^{+}}
$$
holds. Then
$$
g(x) = \dfrac{3 x^5 - 40 \sin x \cos x + 60 \, x \cos x - 20 \sin x}{20 \, x^7},
\hspace*{5 mm} x \in \left(0, \dfrac{\pi}{2} \right) .
$$
Note that
$$
\varphi_{p} (x)=\big{(}p-g(x)\big{)} \frac{x^{6}}{\cos x} \,.
$$

\begin{lemma}
The function $g(x)$ is strictly decreasing for $x\in (0, \pi/2)$.
\end{lemma}
\begin{proof}
Let us notice that the derivative $g'$ is
$$
\!\!
\begin{array}{rcl}
g'(x)
\!\!\!&\!\!=\!\!&\!\!\!
\dfrac{
-40x\cos^2\!x
-
190  x \cos x
+
140 \sin x \cos x
-
\left( 30 x^2 - 70 \right) \sin x
-
3 x^5
+
20 x}{10 x^8}                                                                                       .
\end{array}
$$
It holds
$$
g'(x)<0, \,\,\,x\in \left(0, \dfrac{\pi}{2} \right) \Longleftrightarrow f(x)>0, \,\,\,x\in \left(0, \dfrac{\pi}{2} \right) ,
$$
where
$$
f(x)
=
40\,x\cos^2\!x
+
190 \, x \cos x
-
140 \sin x \cos x
+
\left( 30 \, x^2 - 70 \right) \sin x
+
3 x^5
-
20 \, x.
$$
According to \cite{Malesevic_Makragic_2016}, there exists proof that the MTP function $f$ is positive for $x\!\in \!(0,\pi/2)$. The proof is given in {\bf Appendix A6}.
\end{proof}

\begin{statement}
Let $$
A = \dfrac{3 \pi^5 - 640}{5 \pi^7} = 0.018412 \ldots
\quad\mbox{and}\quad
B = \dfrac{3}{140} = 0.021428 \ldots \,.
$$
Then, it holds$:$

\smallskip
\noindent
{\boldmath $(i)$}
If $p \in (0,A]$, then
$$
\begin{array}{rcl}
x \in \left(0, \dfrac{\pi}{2} \right)
& \Longrightarrow &
2 \, \dfrac{\sin x}{x} + \dfrac{\tan x}{x}
<
3 + \dfrac{3}{20} \, x^4 \! \dfrac{1}{\cos x} - A \, x^6 \dfrac{1}{\cos x} \leq                 \\[2.0 ex]
&                 &
\leq
3 + \dfrac{3}{20} \, x^4 \! \dfrac{1}{\cos x} - p \, x^6 \dfrac{1}{\cos x} \,.
\end{array}
$$

\smallskip
\noindent
{\boldmath $(ii)$}
If $p \in (A,B)$, then $\varphi_{p}(x)$ has exactly one zero $x_{0}^{(p)}$ on $(0,\pi/2)$. Also,
$$
x \in \left(0,x_{0}^{(p)} \right)
\Longrightarrow
2 \, \dfrac{\sin x}{x} + \dfrac{\tan x}{x}
<
3 + \dfrac{3}{20} \, x^4 \! \dfrac{1}{\cos x} - p \, x^6 \dfrac{1}{\cos x}
$$
and
$$
x \in \left(x_{0}^{(p)}, \dfrac{\pi}{2} \right)
\Longrightarrow
3 + \dfrac{3}{20} \, x^4 \! \dfrac{1}{\cos x} - p \, x^6 \dfrac{1}{\cos x}
<
2 \, \dfrac{\sin x}{x} + \dfrac{\tan x}{x} \,
$$
hold.

\smallskip
\noindent
{\boldmath $(iii)$}
If $p \in [B,\infty)$, then
$$
\begin{array}{rcl}
x \in \left(0, \dfrac{\pi}{2} \right)
& \Longrightarrow &
3 + \dfrac{3}{20} \, x^4 \! \dfrac{1}{\cos x} - p \, x^6 \dfrac{1}{\cos x} \leq 3 + \dfrac{3}{20} \, x^4 \! \dfrac{1}{\cos x} - B \, x^6 \dfrac{1}{\cos x}                \\[2.0 ex]
&                 &
<
2 \, \dfrac{\sin x}{x} + \dfrac{\tan x}{x}
\,.
\end{array}
$$
\end{statement}
\begin{proof}
This statement is based on the results of the paper \cite{Mortici_2011} and the fact that
$$
A \;= \lim\limits_{x \rightarrow \pi/2-}{g(x)}
\quad\mbox{and}\quad
B \;=\!\! \lim\limits_{x \rightarrow 0+}{g(x)}\,.
$$
The function $g(x)$ is continuous and, according to Lemma 7, strictly decreasing on $(0,\pi/2)$.

\smallskip
\noindent
{\boldmath $(i)$}
If $p\in (0,A]$, then
$$
g(x)>A \Longleftrightarrow (A-g(x)) \frac{x^{6}}{\cos x}<0 \Longleftrightarrow \varphi_{A}(x)<0
$$
and, therefore, we can conclude that
$$
\begin{array}{rcl}
x \in \left(0, \dfrac{\pi}{2} \right)
& \Longrightarrow &
2 \, \dfrac{\sin x}{x} + \dfrac{\tan x}{x}
<
3 + \dfrac{3}{20} \, x^4 \! \dfrac{1}{\cos x} - A \, x^6 \dfrac{1}{\cos x} \leq                 \\[2.0 ex]
&                 &
\leq
3 + \dfrac{3}{20} \, x^4 \! \dfrac{1}{\cos x} - p \, x^6 \dfrac{1}{\cos x} \,.
\end{array}
$$

\smallskip
\noindent
{\boldmath $(ii)$}
If $p\in (A,B)$, based on Lemma 7, the equation
$$
g(x)=p
$$
has a unique solution $x_{0}^{(p)}$ and it holds
$$
\begin{array}{rcl}
x \in \left(0,x_{0}^{(p)} \right)
\Longrightarrow
g(x)>p
\!\!&\!\!\Longleftrightarrow\!\!&\!\!
\varphi_{p}(x)<0                             \\[3.0 ex]
\!\!&\!\! \Longleftrightarrow\!\!&\!\!
2 \, \dfrac{\sin x}{x} + \dfrac{\tan x}{x}
<
3 + \dfrac{3}{20} \, x^4 \! \dfrac{1}{\cos x} - p \, x^6 \dfrac{1}{\cos x}
\end{array}
$$
and
$$
\begin{array}{rcl}
x \in \left(x_{0}^{(p)}, \dfrac{\pi}{2} \right)
\Longrightarrow
g(x)<p
\!\!&\!\!\Longleftrightarrow\!\!&\!\!
\varphi_{p}(x)>0                             \\[3.0 ex]
\!\!&\!\! \Longleftrightarrow\!\!&\!\!
3 + \dfrac{3}{20} \, x^4 \! \dfrac{1}{\cos x} - p \, x^6 \dfrac{1}{\cos x}
<
2 \, \dfrac{\sin x}{x} + \dfrac{\tan x}{x}.
\end{array}
$$

\smallskip
\noindent
{\boldmath $(iii)$}
If $p\in [B,\infty)$, then
$$
g(x)<B \Longleftrightarrow (B-g(x)) \frac{x^{6}}{\cos x}>0 \Longleftrightarrow \varphi_{B}(x)>0
$$
and, therefore, we can conclude that
$$
\begin{array}{rcl}
x \in \left(0, \dfrac{\pi}{2} \right)
& \Longrightarrow &
3 + \dfrac{3}{20} \, x^4 \! \dfrac{1}{\cos x} - p \, x^6 \dfrac{1}{\cos x} \leq 3 + \dfrac{3}{20} \, x^4 \! \dfrac{1}{\cos x} - B \, x^6 \dfrac{1}{\cos x}                \\[2.0 ex]
&                 &
<
2 \, \dfrac{\sin x}{x} + \dfrac{\tan x}{x}
\,.
\end{array}
$$ 
\end{proof}

Let us notice that $\varphi_{B}(\pi/2-)=+\infty$. Hence, one of the conditions of Theorem 7 is not satisfied, thus, the minimax approximant is not considered.

Based on the previous considerations, enhancement of Theorem 3 has been obtained in the following form:

\begin{proposition}
For every $0 < x < \pi/2$, the following inequalities hold
$$
3 + \left(\dfrac{3 x^4}{20} - \dfrac{3 x^6}{140}\right)\dfrac{1}{\cos x}
\, < \,
2 \,\dfrac{\sin x}{x} + \dfrac{\tan x}{x}
\, < \,
3 + \left(\dfrac{3 x^4}{20} - \dfrac{(3 \pi^5 - 640)x^6}{5 \pi^7}\right)\dfrac{1}{\cos x},
$$
and the constants $A = \dfrac{3 \pi^5 - 640}{5 \pi^7} = 0.018412 \ldots$ and
$B = \dfrac{3}{140} = 0.021428 \ldots \, $ are the best possible.
\end{proposition}

% Statement 4 Statement 4 Statement 4 Statement 4 Statement 4 Statement 4 Statement 4 Statement 4 Statement 4 Statement 4
% ***********************************************************************************************************************

\smallskip\noindent
\underline{\textbf{Improvement of Theorem 4}}
\begin{lemma}
The family of functions
$$
\hspace*{5.00 mm}
\varphi_p(x)
=
\left(\dfrac{\sin x}{x}\right)^2
+ \dfrac{\tan x}{x} - 2 - \dfrac{8}{45} \, x^4 \! \dfrac{1}{\cos x} + p \, x^6 \dfrac{1}{\cos x}
\hspace*{2.25 mm} {\Big (}\mbox{for $x \in \left(0, \dfrac{\pi}{2} \right)$}{\Big )}
$$
is increasingly stratified with respect to parameter $p \!\in\! \RR^{+}$.
\end{lemma}

Let us introduce the function $g(x)$ so that the equivalence
$$
\varphi_{p} (x)=0 \Longleftrightarrow p=g(x), \,\,\,\,\, x\in \left(0, \dfrac{\pi}{2} \right), \,p\in \mathbb{R^{+}}
$$
holds. Then
$$
g(x) = \dfrac{8 x^6 + 90 \, x^2 \cos x - 45 \sin^2 x \, \cos x - 45 \, x \sin x }{45 \, x^8},
\hspace*{5 mm} x \in \left(0, \dfrac{\pi}{2} \right) .
$$
Note that
$$
\varphi_{p} (x)=\big{(}p-g(x)\big{)} \frac{x^{6}}{\cos x} \,.
$$

\begin{lemma}
The function $g(x)$ is strictly decreasing for $x\in (0, \pi/2)$.
\end{lemma}
\begin{proof}
Let us notice that the derivative $g'$ is
$$
\hspace*{-1 mm}
\begin{array}{rcl}
g'(x)
\!\!\!&\!\!\!=\!\!&\!\!\!\!
\dfrac{
-360\cos^3\!x
\!-\!\!
135 x \sin x \cos^2\!x
\!-\!\!\left( 585 x^2\!
-\!360 \right)\! \cos x
\!-\!\!\left( 90 x^3
\!-\!360 x \right) \! \sin x
\!-\!\!16 x^6}{45 x^9} \!.
\end{array}
$$
It holds
$$
g'(x)<0, \,\,\,x\in \left(0, \dfrac{\pi}{2} \right) \Longleftrightarrow f(x)>0, \,\,\,x\in \left(0, \dfrac{\pi}{2} \right) ,
$$
where
$$
f(x)
\!=\!
360 \cos^3\!x
+
135 x \sin x \cos^2\!x
+\left( 585 x^2
-360\right) \cos x
+\left( 90 x^3
-360 x \right) \sin x
+16 x^6.
$$
According to \cite{Malesevic_Makragic_2016}, there exists proof that the MTP function $f$ is positive for $x\!\in \!(0,\pi/2)$. The proof is given in {\bf Appendix A7}.
\end{proof}

\begin{statement}
Let $$
A = \dfrac{32 \pi^5 - 5760}{45 \pi^7} = 0.029670 \ldots
\quad\mbox{and}\quad
B = \dfrac{4}{105} = 0.038095 \ldots \,.
$$
Then, it holds$:$

\smallskip
\noindent
{\boldmath $(i)$}
If $p \in (0,A]$, then
$$
\begin{array}{rcl}
x \in \left(0, \dfrac{\pi}{2} \right)
& \Longrightarrow &
\left(\dfrac{\sin x}{x}\right)^2
+ \dfrac{\tan x}{x}
<
2 + \dfrac{8}{45} \, x^4 \! \dfrac{1}{\cos x} - A \, x^6 \dfrac{1}{\cos x} \leq                 \\[2.0 ex]
&                 &
\leq
2 + \dfrac{8}{45} \, x^4 \! \dfrac{1}{\cos x} - p \, x^6 \dfrac{1}{\cos x} \,.
\end{array}
$$

\smallskip
\noindent
{\boldmath $(ii)$}
If $p \in (A,B)$, then $\varphi_{p}(x)$ has exactly one zero $x_{0}^{(p)}$ on $(0, \pi/2)$. Also,
$$
x \in \left(0,x_{0}^{(p)} \right)
\Longrightarrow
\Big{(}\dfrac{\sin x}{x}\Big{)}^{2} + \dfrac{\tan x}{x}
<
2 + \dfrac{8}{45} \, x^4 \! \dfrac{1}{\cos x} - p \, x^6 \dfrac{1}{\cos x}
$$
and
$$
x \in \left(x_{0}^{(p)}, \dfrac{\pi}{2} \right)
\Longrightarrow
2 + \dfrac{8}{45} \, x^4 \! \dfrac{1}{\cos x} - p \, x^6 \dfrac{1}{\cos x}
<
\Big{(}\dfrac{\sin x}{x}\Big{)}^{2} + \dfrac{\tan x}{x}
$$
hold.

\smallskip
\noindent
{\boldmath $(iii)$}
If $p \in [B,\infty)$, then
$$
\begin{array}{rcl}
x \in \left(0, \dfrac{\pi}{2} \right)
& \Longrightarrow &
2 + \dfrac{8}{45} \, x^4 \! \dfrac{1}{\cos x} - p \, x^6 \dfrac{1}{\cos x} \leq 2 + \dfrac{8}{45} \, x^4 \! \dfrac{1}{\cos x} - B \, x^6 \dfrac{1}{\cos x}                \\[2.0 ex]
&                 &
<
\Big{(}\dfrac{\sin x}{x}\Big{)}^{2} + \dfrac{\tan x}{x}
\,.
\end{array}
$$
\end{statement}
\begin{proof}
This statement is based on the results of the paper \cite{Mortici_2011} and the fact that
$$
A \;= \lim\limits_{x \rightarrow \pi/2-}{g(x)}
\quad\mbox{and}\quad
B \;=\!\! \lim\limits_{x \rightarrow 0+}{g(x)} \,.
$$

\noindent
The function $g(x)$ is continuous and, according to Lemma 9, strictly decreasing on $(0,\pi/2)$.

\smallskip
\noindent
{\boldmath $(i)$}
If $p\in (0,A]$, then
$$
g(x)>A \Longleftrightarrow (A-g(x)) \frac{x^{6}}{\cos x}<0 \Longleftrightarrow \varphi_{A}(x)<0
$$
and, therefore, we can conclude that
$$
\begin{array}{rcl}
x \in \left(0, \dfrac{\pi}{2} \right)
& \Longrightarrow &
\left(\dfrac{\sin x}{x}\right)^2
+ \dfrac{\tan x}{x}
<
2 + \dfrac{8}{45} \, x^4 \! \dfrac{1}{\cos x} - A \, x^6 \dfrac{1}{\cos x} \leq                 \\[2.0 ex]
&                 &
\leq
2 + \dfrac{8}{45} \, x^4 \! \dfrac{1}{\cos x} - p \, x^6 \dfrac{1}{\cos x} \,.
\end{array}
$$

\smallskip
\noindent
{\boldmath $(ii)$}
If $p\in (A,B)$, based on Lemma 9, the equation
$$
g(x)=p
$$
has a unique solution $x_{0}^{(p)}$ and it holds
$$
\begin{array}{rcl}
x \in \left(0,x_{0}^{(p)} \right)
\Longrightarrow
g(x)>p
\!\!&\!\!\Longleftrightarrow\!\!&\!\!
\varphi_{p}(x)<0                             \\[3.0 ex]
\!\!&\!\! \Longleftrightarrow\!\!&\!\!
\Big{(}\dfrac{\sin x}{x}\Big{)}^{2} + \dfrac{\tan x}{x}
<
2 + \dfrac{8}{45} \, x^4 \! \dfrac{1}{\cos x} - p \, x^6 \dfrac{1}{\cos x}
\end{array}
$$
and
$$
\begin{array}{rcl}
x \in \left(x_{0}^{(p)}, \dfrac{\pi}{2} \right)
\Longrightarrow
g(x)<p
\!\!&\!\!\Longleftrightarrow\!\!&\!\!
\varphi_{p}(x)>0                             \\[3.0 ex]
\!\!&\!\! \Longleftrightarrow\!\!&\!\!
2 + \dfrac{8}{45} \, x^4 \! \dfrac{1}{\cos x} - p \, x^6 \dfrac{1}{\cos x}
<
\Big{(}\dfrac{\sin x}{x}\Big{)}^{2} + \dfrac{\tan x}{x}.
\end{array}
$$

\smallskip
\noindent
{\boldmath $(iii)$}
If $p\in [B,\infty)$, then
$$
g(x)<B \Longleftrightarrow (B-g(x))= \frac{x^{6}}{\cos x}>0 \Longleftrightarrow \varphi_{B}(x)>0
$$
and, therefore, we can conclude that
$$
\begin{array}{rcl}
x \in \left(0, \dfrac{\pi}{2} \right)
& \Longrightarrow &
2 + \dfrac{8}{45} \, x^4 \! \dfrac{1}{\cos x} - p \, x^6 \dfrac{1}{\cos x} \leq 2 + \dfrac{8}{45} \, x^4 \! \dfrac{1}{\cos x} - B \, x^6 \dfrac{1}{\cos x}                \\[2.0 ex]
&                 &
<
\Big{(}\dfrac{\sin x}{x}\Big{)}^{2} + \dfrac{\tan x}{x}
\,.
\end{array}
$$ 
\end{proof}

Let us notice that $\varphi_{B}(\pi/2-)=+\infty$. Hence, one of the conditions of Theorem 7 is not satisfied, thus, the minimax approximant is not considered.

Based on the previous considerations, enhancement of Theorem 4 has been obtained in the following form:

\begin{proposition}
For every $0 < x < \pi/2$, the following inequalities hold
$$
2 + \left(\dfrac{8 x^4}{45} - \dfrac{4 x^6}{105}\right)\dfrac{1}{\cos x}
<
\left(\dfrac{\sin x}{x}\right)^2
+ \dfrac{\tan x}{x}
<
2 + \left(\dfrac{8 x^4}{45} - \dfrac{(32 \pi^5 - 5760)x^6}{45 \pi^7}\right)\dfrac{1}{\cos x},
$$
and the constants $A = \dfrac{32 \pi^5 - 5760}{45 \pi^7} = 0.029670 \ldots$ and
$B = \dfrac{4}{105} = 0.038095 \ldots \, $ are the best possible.
\end{proposition}

\break

Analogously to Statement 1 and Statement 2, the following statements can be proved:

% Statement 5 Statement 5 Statement 5 Statement 5 Statement 5 Statement 5 Statement 5 Statement 5 Statement 5 Statement 5
% ***********************************************************************************************************************

\bigskip\noindent
\underline{\textbf{Improvement of Theorem 5}}
\begin{lemma}
The family of functions
$$
\hspace*{35.0 mm}
\varphi_p(x)
=
\left(\dfrac{x}{\sin x}\right)^2
\, + \, \dfrac{x}{\tan x} - 2 - \, p \, x^4
\hspace*{8.80 mm} {\Big (}\mbox{for $x \in \left(0, \dfrac{\pi}{2} \right)$}{\Big )}
$$
is decreasingly stratified with respect to parameter $p \!\in\! \RR^{+}$.
\end{lemma}

\begin{statement}
Let $$
A = \dfrac{2}{45} = 0.044444 \ldots
\quad\mbox{and}\quad
B = \dfrac{4 \pi^2-32}{\pi^4} = 0.076773 \ldots \,.
$$
Then, it holds$:$

\smallskip
\noindent
{\boldmath $(i)$}
If $p \in (0,A]$, then
$$
x \in \left(0, \dfrac{\pi}{2} \right)
\Longrightarrow
2 + p x^4
\leq
2 + A x^4
<
\left(\dfrac{x}{\sin x}\right)^2
\, + \, \dfrac{x}{\tan x}.
$$

\smallskip
\noindent
{\boldmath $(ii)$}
If $p \in (A,B)$, then $\varphi_{p}(x)$ has exactly one zero $x_{0}^{(p)}$ on $(0,\pi/2)$. Also,
$$
x \in \left(0,x_{0}^{(p)} \right)
\Longrightarrow
\left(\dfrac{x}{\sin x}\right)^2
\, + \, \dfrac{x}{\tan x}
<
2 + p \, x^4
$$
and
$$
x \in \left(x_{0}^{(p)}, \dfrac{\pi}{2} \right)
\Longrightarrow
2 + p \, x^4
<
\left(\dfrac{x}{\sin x}\right)^2
\, + \, \dfrac{x}{\tan x}
$$
hold.

\smallskip
\noindent
{\boldmath $(iii)$}
If $p \in [B,\infty)$, then
$$
x \in \left(0, \dfrac{\pi}{2} \right)
\Longrightarrow
\left(\dfrac{x}{\sin x}\right)^2
\, + \, \dfrac{x}{\tan x}
<
2 + B x^4
\leq
2 + p x^4.
$$

\smallskip
\noindent
{\boldmath $(iv)$}
There is exactly one solution to the following equation
$$
\left| \varphi_p{\big (}t^{(p)}{\big )}\right| =  \varphi_p \left(\frac{\pi}{2}- \right) ,
$$
where $t^{(p)}$ is a minimum of  $\varphi_p(x)$ on $(0, \pi/2)$, with respect to parameter $p \!\in\! (A,B)$, which is numerically determined as
$$
p_0 = 0.072425 \ldots \,.
$$
For the value
$$
d_{0}
=
\left| \varphi_{p_0}{\big (}t^{(p_0)}{\big )} \right|=  \varphi_{p_0} \left(\frac{\pi}{2}- \right)
=
0.026471\ldots \,,
$$
the following holds
$$
d_0 = \inf\limits_{p \in \RR^{+}} \sup\limits_{x \in (0,\pi/2)}{\left| \varphi_p(x) \right|}.
$$

\smallskip
\noindent
{\boldmath $(v)$}
For the value $p_0 = 0.072425 \ldots  $ , the minimax approximant of the family is determined as
$$
\varphi_{p_0}(x)
=
\left(\dfrac{x}{\sin x}\right)^2
\, + \, \dfrac{x}{\tan x}
-
2
-
0.072425 \ldots  x^4
$$
and it determines the corresponding minimax approximation
\begin{equation*}
\label{minimax_approximation_T5}
\left(\dfrac{x}{\sin x}\right)^2
\, + \, \dfrac{x}{\tan x}
\approx
2 + 0.072425\ldots x^4 \,.
\end{equation*}
\end{statement}

Based on the previous considerations, enhancement of Theorem 5 has been obtained in the following form:

\begin{proposition}
For every $0 < x < \pi/2$, the following inequalities hold
$$
2 + \dfrac{2}{45}x^4
\,<\,
\left(\dfrac{x}{\sin x}\right)^2
\, + \, \dfrac{x}{\tan x}
\,<\,
2+ \dfrac{4 \pi^2-32}{\pi^4} \,x^4,
$$
and the constants $A = \dfrac{2}{45} = 0.044444 \ldots $ and $B = \dfrac{4 \pi^2-32}{\pi^4} = 0.076773 \ldots$
are the best possible.
\end{proposition}

% Statement 6 Statement 6 Statement 6 Statement 6 Statement 6 Statement 6 Statement 6 Statement 6 Statement 6 Statement 6
% ***********************************************************************************************************************

\smallskip\noindent
\underline{\textbf{Improvement of Theorem 6}}
\begin{lemma}
The family of functions
$$
\hspace*{25.0 mm}
\varphi_p(x)
=
3 \, \dfrac{x}{\sin x}
\, + \, \cos x \, - \, 4 \, - \, \dfrac{1}{10} \, x^4 \, -  \, p \, x^6
\hspace*{8.80 mm} {\Big (}\mbox{for $x \in \left(0, \dfrac{\pi}{2} \right)$}{\Big )}
$$
is decreasingly stratified with respect to parameter $p \!\in\! \RR^{+}$.
\end{lemma}

\begin{statement}
Let $$
A = \dfrac{1}{210} = 0.0047619 \ldots
\quad\mbox{and}\quad
B = \dfrac{480 \pi \!-\! 2 \pi^4 \!-\! 1280}{5\pi^6} = 0.0068954 \ldots \,.
$$
Then, it holds$:$

\smallskip
\noindent
{\boldmath $(i)$}
If $p \in (0,A]$, then
$$
x \in \left(0, \dfrac{\pi}{2} \right)
\Longrightarrow
4 + \dfrac{1}{10} x^4 + p x^6
\leq
4 + \dfrac{1}{10} x^4 + A x^6
<
3 \, \dfrac{x}{\sin x}
\, + \, \cos x.
$$

\smallskip
\noindent
{\boldmath $(ii)$}
If $p \in (A,B)$, then $\varphi_{p}(x)$ has exactly one zero $x_{0}^{(p)}$ on $(0,\pi/2)$. Also,
$$
x \in \left(0,x_{0}^{(p)} \right)
\Longrightarrow
3 \, \dfrac{x}{\sin x}
\, + \, \cos x
<
4 + \dfrac{1}{10} x^4 + p x^6
$$
and
$$
x \in \left(x_{0}^{(p)}, \dfrac{\pi}{2} \right)
\Longrightarrow
4 + \dfrac{1}{10} x^4 + p x^6
<
3 \, \dfrac{x}{\sin x}
\, + \, \cos x
$$
hold.

\smallskip
\noindent
{\boldmath $(iii)$}
If $p \in [B,\infty)$, then
$$
x \in \left(0, \dfrac{\pi}{2} \right)
\Longrightarrow
3 \, \dfrac{x}{\sin x}
\, + \, \cos x
<
4 + \dfrac{1}{10} x^4 + B x^6
\leq
4 + \dfrac{1}{10} x^4 + p x^6.
$$

\smallskip
\noindent
{\boldmath $(iv)$}
There is exactly one solution to the following equation
$$
\left| \varphi_p{\big (}t^{(p)}{\big )}\right| =  \varphi_p \left(\frac{\pi}{2}- \right) ,
$$
where $t^{(p)}$ is a minimum of  $\varphi_p(x)$ on $(0, \pi/2)$, with respect to parameter $p \!\in\! (A,B)$, which is numerically determined as
$$
p_0 = 0.0066982 \ldots \,.
$$
For the value
$$
d_{0}
=
\left| \varphi_{p_0}{\big (}t^{(p_0)}{\big )} \right|=  \varphi_{p_0} \left(\frac{\pi}{2}- \right)
=
0.0029637 \ldots \,,
$$
the following holds
$$
d_0 = \inf\limits_{p \in \RR^{+}} \sup\limits_{x \in (0,\pi/2)}{\left| \varphi_p(x) \right|}.
$$

\smallskip
\noindent
{\boldmath $(v)$}
For the value $p_0 = 0.0066982 \ldots  $ , the minimax approximant of the family is determined as
$$
\varphi_{p_0}(x)
=
3 \, \dfrac{x}{\sin x}
+\cos x
-
4
-
\dfrac{1}{10} \, x^4
-
0.0066982 \ldots  x^6
$$
and it determines the corresponding minimax approximation
\begin{equation*}
\label{minimax_approximation_T6}
3 \, \dfrac{x}{\sin x}
+\cos x
\approx
4+\dfrac{1}{10} \, x^4
+0.0066982 \ldots  x^6.
\end{equation*}
\end{statement}

Based on the previous considerations, enhancement of Theorem 6 has been obtained in the following form:

\begin{proposition}
For every $0 < x < \pi/2$, the following inequalities hold
$$
4 + \dfrac{1}{10} x^4 + \dfrac{1}{210} x^6
\,<\,
3 \, \dfrac{x}{\sin x}
\, + \, \cos x
\,<\,
4 + \dfrac{1}{10} x^4 + \dfrac{480 \pi - 2 \pi^4 - 1280}{5\pi^6} \,x^6,
$$
and the constants $A = \dfrac{1}{210} = 0.0047619 \ldots $ and $B = \dfrac{480 \pi - 2 \pi^4 - 1280}{5\pi^6} = 0.0068954 \ldots$
are the best possible.
\end{proposition}

The existence of minimax approximant in Statement 5 and Statement 6 is a consequence of Theorem 7' and Theorem 8.

\break

\vspace{1.5cc}
\begin{section}
{CONCLUSION}
\end{section}

This paper specifies the results of C. Mortici \cite{Mortici_2011} using the method described in \cite{Malesevic_Mihailovic_2021}. The mentioned method could be applied for possible improvements of existing results from the Theory of analytic inequalities \cite{Mitrinovic_1970}, \cite{Milovanovic_Rassias_2014}, \cite{Cloud_Drachman_Lebedev_2014} in terms of determining the corresponding minimax approximants for various inequalities. In the previous section, examples of minimax approximations are presented where they exist.
It is important to note that through minimax approximants, the error in approximations is minimized in the considered sense.

The main aim of this paper is to promote SimTheP, an automated theorem prover for MTP inequalities, developed through the doctoral dissertation \cite{Banjac_2019}.
All the essential proofs of MTP inequalities in this paper are given in the Appendix and derived using the prover SimTheP.
For any given MTP inequality $f(x) > 0$, for $x \!\in\! \SSS \subseteq [0, \pi/2]$, SimTheP provides a structured proof divided into parts \textbf{I}-\textbf{IV}.
Each part is designed to allow manual step-by-step verification, demonstrating SimTheP's capability to replicate the human way of proving inequalities.

Many results within the Theory of analytic inequalities, for example, results
from \cite{Malesevic_Mihailovic_2021}, \cite{Malesevic_Makragic_2016}, \cite{Micovic_Malesevic_2024}-\cite{Malesevic_Mihailovic_NenezicJovic_Milinkovic_2022}, \cite{Shiping_Zhong_2016}-\cite{Milovanovic_Rassias_2014}, \cite{Banjac_Makragic_Malesevic_2016}, \cite{Alzer_Kwong_2016}, \cite{Alzer_Kwong_2020}-\cite{Huy_Hieu_Van_2024}, could be proved
using the prover SimTheP, see the link https://simthep.etf.bg.ac.rs/.

It is crucial to highlight that through Statements 1--6, all Theorems 1--6 have been improved and minimax approximations have been determined wherever feasible. As a result, Propositions 1--6 were obtained, where for the inequalities considered in Theorems 1--6, the best possible constants were identified. Such an approach to Theorems 1--6 was enabled by the utilization of the concept of stratification.
Moreover, this paper presents the first paper in which the automated theorem prover SimTheP was utilized to deliver proofs for the MTP inequalities, marking a significant advancement in the field.

\bigskip
\noindent
{\bf Acknowledgments.}

The authors are greatly indebted to Dr. Ivana Jovovi\' c for numerous stimulating conversations about the concept of stratification. Her assistance in applying MTP inequalities to various problems, coupled with her ongoing and steadfast commitment to supporting our endeavors, has been instrumental and is deeply appreciated. Special appreciation is also extended to Dr. Marija Nenezi\' c Jovi\' c for her invaluable comments that greatly enhanced certain proofs. The authors also express their gratitude to the referees for their thorough reading of the paper and their valuable suggestions and comments.

This work was financially supported by the Ministry of Science, Technological Development
and Innovation of the Republic of Serbia under contract numbers:
451-03-65/2024-03/200156 (for the first author),
451-03-65/2024-03/200103 (for the second, fourth and fifth authors)
and 451-03-66/2024-03/200103 (for the third author).
The research of the first author has also been supported by the
Faculty of Technical Sciences, University of Novi Sad through project
"Improving the teaching process in the English language in fundamental
disciplines" (No. 01-3394/1).

\vspace{1cc}

\newpage

\bigskip\medskip
\noindent
\centerline{\small \bf APPENDIX}

\bigskip
{\small

This Appendix was created using the automated prover SimTheP, which for the MTP function and the interval gives as output TeX/PDF files that were directly transferred to the Appendix.

\bigskip\noindent
{\bf APPENDIX A1} % A1 A1 A1 A1 A1 A1 A1 A1 A1 A1 A1 A1 A1 A1 A1 A1 A1 A1 A1 A1 A1 A1 A1 A1 A1 A1 A1 A1 A1

\medskip\noindent
The initial MTP function is
$$
\begin{array}{rcl}
f(x)
\!\!&\!\!=\!\!&\!\!
2\cdot x^{7}+135\cdot\sin x \cdot\cos^{2}x+\left(15\cdot x^{4} -135\right)\cdot\sin x
\\[1.25ex]
\!\!&\!\! \!\!&\!\!
-45\cdot x\cdot\cos^{3}x +\left(90\cdot x^{3} +45\cdot x\right)\cdot\cos x
\end{array}
$$
and the initial interval is $\SSS = \left(0, \frac{\pi}{2} \right)$.

\medskip\noindent
\underline{\em Automated proof that $f(x) > 0$ for $x \!\in\! \left(0, \frac{\pi}{2} \right)$$:$}

\medskip
\noindent
\textbf{I} (Recognition of possible case)
Facts
$f\left(0\right) = 0$
and
$f\left(\frac{\pi}{2}\right) \neq 0$
are correct.
The MTP function $f$ is positive at boundary point $\frac{\pi}{2}$ {\big (} $\!f\left(\frac{\pi}{2}\right) = 3.51310 \ldots  >0\,${\big )}.
Therefore, it is possible that $f(x) > 0$ over $\left(0,\frac{\pi}{2}\right)$.

\medskip
\noindent
\textbf{II} (Transformation of angles)
After the transformation of terms $\cos^m x \cdot \sin^n x$ $(m, n \!\in\! N_0)$
into the sum of sine and cosine functions of multiple angles, in the MTP function $f$,
we obtain
$$
\begin{array}{rcl}
f(x)
\!\!&\!\!=\!\!&\!\!
-\frac{45}{4}\cdot x\cdot\cos 3 x +\left(15\cdot x^{4} -\frac{405}{4}\right)\cdot\sin x
\\[1.25ex]
\!\!&\!\! \!\!&\!\!
+ \left(90\cdot x^{3}+\frac{45}{4}\cdot x\right)\cdot\cos x +2\cdot x^{7}+\frac{135}{4}\cdot\sin 3x \,.
\end{array}
$$
Then, we consider the previous expression as two separate expressions, the first with positive and the second with negative terms next to sine and cosine functions
$$
\begin{array}{l}
\begin{array}{rcl}
f^{+}(x)
\!\!&\!\!=\!\!&\!\!
\frac{135}{4}\cdot\sin 3 x +\left(90\cdot x^{3} +\frac{45}{4}\cdot x\right)\cdot\cos x
+15\cdot x^{4}\cdot\sin x
+2\cdot x^{7} \,,
\end{array} \\[2.50 ex]
\begin{array}{rcl}
f^{-}(x)
\!\!&\!\!=\!\!&\!\!
-\frac{405}{4}\cdot\sin x -\frac{45}{4}\cdot x\cdot\cos 3x \,.
\end{array}
\end{array}
$$

\medskip\noindent
\textbf{III} (Determination of downward rational polynomial approximation)
After substitution of sine and cosine functions by appropriate (downward or upward)
polynomial approximations, we obtain downward polynomial approximations of $f^{+}(x)$ and $f^{-}(x)$
respectively
$$
\begin{array}{l}
\begin{array}{rcl}
P^{+}_{i_{0},i_{1},i_{2}}(x)
\!\!&\!\!=\!\!&\!\!
\frac{135}{4}\cdot \underline{T}_{\,4\cdot i_{0}+3}^{\,\sin,0}(3 x)+\left(90\cdot x^{3} +\frac{45}{4}\cdot x\right)\cdot \underline{T}_{\,4\cdot i_{1}+2}^{\,\cos,0}(x)   \\[1.25ex]
\!\!&\!\! \!\!&\!\!
+15\cdot x^{4}\cdot \underline{T}_{\,4\cdot i_{2}+3}^{\,\sin,0}(x)+2\cdot x^{7} \,,
\\[2.25ex]
\end{array} \\[3.00 ex]
\begin{array}{rcl}
P^{-}_{i_{3},i_{4}}(x)
\!\!&\!\!=\!\!&\!\!
-\frac{405}{4}\cdot \overline{T}_{\,4\cdot i_{3}+1}^{\,\sin,0}(x)-\frac{45}{4}\cdot x\cdot \overline{T}_{\,4\cdot i_{4}}^{\,\cos,0}(3 x)\,.
\end{array}
\end{array}
$$
For concrete indices $\left(i_{0},i_{1},i_{2},i_{3},i_{4} \right)=\left(2, 2, 1, 2, 3\right)$, we obtain
$$
\begin{array}{l}
\begin{array}{rcl}
P^{+}_{2,2,1}(x)
\!\!&\!\!=\!\!&\!\!
\frac{135}{4}\cdot \underline{T}_{\,4\cdot 2+3}^{\,\sin,0}(3 x)
+
\left(90\cdot x^{3} +\frac{45}{4}\cdot x\right)\cdot \underline{T}_{\,4\cdot 2+2}^{\,\cos,0}(x)  \\[1.25ex]
\!\!&\!\! \!\!&\!\!
+15\cdot x^{4}\cdot \underline{T}_{\,4\cdot 1+3}^{\,\sin,0}(x)+2\cdot x^{7} \,,
\\[2.25ex]
\end{array} \\[3.00 ex]
\begin{array}{rcl}
P^{+}_{2,2,1}(x)
\!\!&\!\!=\!\!&\!\!
-\frac{1}{40320}\cdot x^{13}-\frac{133523}{887040}\cdot x^{11}+\frac{3281}{1792}\cdot x^{9} \\[1.25ex]
\!\!&\!\! \!\!&\!\!
-\frac{639}{56}\cdot x^{7} +\frac{621}{16}\cdot x^{5}-\frac{135}{2}\cdot x^{3}+\frac{225}{2}\cdot x
\end{array}
\end{array}
$$
and
$$
\begin{array}{l}
\begin{array}{rcl}
P^{-}_{2,3}(x)
\!\!&\!\!=\!\!&\!\!
-\frac{405}{4}\cdot \overline{T}_{\,4\cdot 2+1}^{\,\sin,0}(x)
-\frac{45}{4}\cdot x\cdot \overline{T}_{\,4\cdot 3}^{\,\cos,0}(3 x) \,,
\end{array} \\[2.50 ex]
\begin{array}{rcl}
P^{-}_{2,3}(x)
\!\!&\!\!=\!\!&\!\!
-\frac{19683}{1576960}\cdot x^{13}+\frac{6561}{35840}\cdot x^{11}-\frac{3281}{1792}\cdot x^{9}\\[1.25 ex]
\!\!&\!\! \!\!&\!\!
+\frac{639}{56}\cdot x^{7}-\frac{621}{16}\cdot x^{5}+\frac{135}{2}\cdot x^{3}-\frac{225}{2}\cdot x \,.
\end{array}
\end{array}
$$
Finally, for the MTP function
$$
f(x)
=
f^{+}(x)+f^{-}(x)
$$
we obtain the concrete {\em downward polynomial approximation}
$$
\begin{array}{rcl}
P(x)
\!\!&\!\!=\!\!&\!\!
P^{+}_{2,2,1}(x)+P^{-}_{2,3}(x) \\[1.5 ex]
\!\!&\!\!=\!\!&\!\!
-\frac{25357}{2027520}\cdot x^{13}+\frac{115447}{3548160}\cdot x^{11}
\end{array}
$$
over $\left(0,\frac{\pi}{2}\right)$, i.e. it holds that
$$
f(x) > P(x)
$$
over $\left(0,\frac{\pi}{2}\right)$.

\medskip
\noindent
\textbf{IV} (The final part)
Based on the Sturm theorem, the following inequality
$$
P(x)>0
$$
is true over $\left(0,\frac{\pi}{2}\right)$.
The stated conclusion for the polynomial function $P$ is correct based on the following facts$:$\\
$1.$
We can conclude, by Sturm theorem, that the polynomial function $P(x)$ has only one zero over the concrete extended segment $[-0.1, 1.58]$
of the initial interval $\left(0,\frac{\pi}{2}\right)$.\\
$2.$
Facts
$P\left(-0.1\right) \neq 0$, $P\left(0\right) = 0$ and $P\left(1.58 \right) \neq 0$
are correct.\\
$3.$ The polynomial $P$ is positive at boundary point $\frac{\pi}{2}$
{\big (}$P\left(\frac{\pi}{2}\right) = 0.24116 \ldots > 0\,${\big )}.
$\!\!$\\[0.1 ex]
Therefore, the following inequality
$$
f(x)>0
$$
is true over  $\left(0,\frac{\pi}{2}\right)$. \hfill $\Box$

\bigskip\noindent
{\bf APPENDIX A2} % A2 A2 A2 A2 A2 A2 A2 A2 A2 A2 A2 A2 A2 A2 A2 A2 A2 A2 A2 A2 A2 A2 A2 A2 A2 A2 A2 A2 A2

\medskip\noindent
The initial MTP function is
$$
\begin{array}{rcl}
f(x)
\!\!&\!\!=\!\!&\!\!
\left(-15309\cdot x^{6} +170100\cdot x^{4} -476280\cdot x^{2} +181440\right)\cdot\sin x \cdot\cos^{2} x                                                     \\[1.25ex]
\!\!&\!\! \!\!&\!\!
+\left(-x^{10} +3843\cdot x^{6} -44100\cdot x^{4} +158760\cdot x^{2} -181440\right)\cdot\sin x                                                         \\[1.25ex]
\!\!&\!\! \!\!&\!\!
+\left(2187\cdot x^{7} -61236\cdot x^{5}+340200\cdot x^{3} -423360\cdot x\right)\cdot\cos^{3} x                                                       \\[1.25ex]
\!\!&\!\! \!\!&\!\!
+\left(-1641\cdot x^{7} +46116\cdot x^{5} -264600\cdot x^{3} +423360\cdot x\right)\cdot\cos x
\end{array}
$$
and the initial interval is $\SSS = \left(0, 1 \right]$.

\medskip\noindent
\underline{\em Automated proof that $f(x) > 0$ for $x \!\in\! \left( 0, 1 \right]$$:$}

\medskip
\noindent
\textbf{I} (Recognition of possible case)
Facts
$f\left(0\right) = 0$
and
$f\left(1\right) \neq 0$
are correct.
The MTP function $f$ is positive at boundary point $1$ {\big (}$ \; f\left(1\right) = 14.68957 \ldots  >0$ {\big )}.
Therefore, it is possible that $f(x) > 0$ over $\left( 0, 1 \right]$.

\medskip
\noindent
\textbf{II} (Transformation of angles)
After the transformation of terms $\cos^m x \cdot \sin^n x$ $(m, n \!\in\! N_0)$
into the sum of sine and cosine functions of multiple angles, in the MTP function $f$,
we obtain
$$
\begin{array}{rcl}
f(x)
\!\!&\!\!=\!\!&\!\!
\left(-\frac{15309}{4}\cdot x^{6} +42525\cdot x^{4} -119070\cdot x^{2} +45360\right)\cdot\sin 3x      \\[1.25ex]
\!\!&\!\! \!\!&\!\!
+\left(\frac{2187}{4}\cdot x^{7} -15309\cdot x^{5} +85050\cdot x^{3} -105840\cdot x\right)\cdot\cos 3x \\[1.25ex]
\!\!&\!\! \!\!&\!\!
+ \left(-x^{10} +\frac{63}{4}\cdot x^{6} -1575\cdot x^{4} +39690\cdot x^{2} -136080\right)\cdot\sin x         \\[1.25ex]
\!\!&\!\! \!\!&\!\!
+\left(-\frac{3}{4}\cdot x^{7} +189\cdot x^{5} -9450\cdot x^{3} +105840\cdot x\right)\cdot\cos x \,.
\end{array}
$$
Then, we consider the previous expression as two separate expressions, the first with positive and the second with negative terms next to sine and cosine functions
$$
\begin{array}{l}
\begin{array}{rcl}
f^{+}(x)
\!\!&\!\!=\!\!&\!\!
\left(42525\cdot x^{4} +45360\right)\cdot\sin 3 x +\left(189\cdot x^{5} +105840\cdot x\right)\cdot\cos x  \\[1.25ex]
\!\!&\!\! \!\!&\!\!
+\left(\frac{63}{4}\cdot x^{6} +39690\cdot x^{2}\right)\cdot\sin x +\left(\frac{2187}{4}\cdot x^{7}
+85050\cdot x^{3}\right)\cdot\cos 3 x \,, \\[2.25ex]
\end{array} \\[3.00 ex]
\begin{array}{rcl}
f^{-}(x)
\!\!&\!\!=\!\!&\!\!
\left(-x^{10} -1575\cdot x^{4} -136080\right)\cdot\sin x
+\left(-15309\cdot x^{5} -105840\cdot x\right)\cdot\cos 3 x                                             \\[1.25ex]
\!\!&\!\! \!\!&\!\!
+\left(-\frac{15309}{4}\cdot x^{6} -119070\cdot x^{2}\right)\cdot\sin 3 x +\left(-\frac{3}{4}\cdot x^{7} -9450\cdot x^{3}\right)\cdot\cos x \,.
\end{array}
\end{array}
$$

\medskip\noindent
\textbf{III} (Determination of downward rational polynomial approximation)
After substitution of sine and cosine functions by appropriate (downward or upward)
polynomial approximations, we obtain downward polynomial approximations of $f^{+}(x)$ and $f^{-}(x)$
respectively
$$
\begin{array}{l}
\begin{array}{rcl}
P^{+}_{i_{0},i_{1},i_{2},i_{3}}(x)
\!\!&\!\!=\!\!&\!\!
\left(42525\cdot x^{4} +45360\right)\cdot \underline{T}_{\,4\cdot i_{0}+3}^{\,\sin,0}(3 x)
 \\[1.25ex]
\!\!&\!\! \!\!&\!\!
+ \left(189\cdot x^{5} +105840\cdot x\right)\cdot \underline{T}_{\,4\cdot i_{1}+2}^{\,\cos,0}(x)  \\[1.25ex]
\!\!&\!\! \!\!&\!\!
+ \left(\frac{63}{4}\cdot x^{6}+39690\cdot x^{2}\right)\cdot \underline{T}_{\,4\cdot i_{2}+3}^{\,\sin,0}(x)  \\[1.25ex]
\!\!&\!\! \!\!&\!\!
+\left(\frac{2187}{4}\cdot x^{7} +85050\cdot x^{3}\right)\cdot \underline{T}_{\,4\cdot i_{3}+2}^{\,\cos,0}(3 x) \,,
\\[2.25ex]
\end{array} \\[3.00 ex]
\begin{array}{rcl}
P^{-}_{i_{4},i_{5},i_{6},i_{7}}(x)
\!\!&\!\!=\!\!&\!\!
\left(-x^{10} -1575\cdot x^{4} -136080\right)\cdot \overline{T}_{\,4\cdot i_{4}+1}^{\,\sin,0}(x)\\[1.25ex]
\!\!&\!\! \!\!&\!\!\
+\left(-15309\cdot x^{5}-105840\cdot x\right)\cdot \overline{T}_{\,4\cdot i_{5}}^{\,\cos,0}(3x)  \\[1.25ex]
\!\!&\!\! \!\!&\!\!\
+\left(-\frac{15309}{4}\cdot x^{6}-119070\cdot x^{2}\right)\cdot \overline{T}_{\,4\cdot i_{6}+1}^{\,\sin,0}(3x)  \\[1.25ex]
\!\!&\!\! \!\!&\!\!\
+\left(-\frac{3}{4}\cdot x^{7}-9450\cdot x^{3}\right)\cdot \overline{T}_{\,4\cdot i_{7}}^{\,\cos,0}(x) \,.
\end{array}
\end{array}
$$
For concrete indices $\left(i_{0},i_{1},i_{2},i_{3},i_{4},i_{5},i_{6},i_{7}\right)=\left(3,3,3,3,3,3,3,3\right)$, we obtain
$$
\begin{array}{l}
\begin{array}{rcl}
P^{+}_{3,3,3,3}(x)
\!\!&\!\!=\!\!&\!\!
\left(42525\cdot x^{4} +45360\right)\cdot \underline{T}_{\,4\cdot 3+3}^{\,\sin,0}(3 x)
 \\[1.25ex]
\!\!&\!\! \!\!&\!\!\
+\left(189\cdot x^{5}+105840\cdot x\right)\cdot \underline{T}_{\,4\cdot 3+2}^{\,\cos,0}(x)
 \\[1.25ex]
\!\!&\!\! \!\!&\!\!
+\left(\frac{63}{4}\cdot x^{6} +39690\cdot x^{2}\right)\cdot \underline{T}_{\,4\cdot 3+3}^{\,\sin,0}(x) \\[1.25ex]
\!\!&\!\! \!\!&\!\!+\left(\frac{2187}{4}\cdot x^{7} +85050\cdot x^{3}\right)\cdot \underline{T}_{\,4\cdot 3+2}^{\,\cos,0}(3 x) \,,
\\[2.25ex]
\end{array} \\[3.00 ex]
\begin{array}{rcl}
P^{+}_{3,3,3,3}(x)
\!\!&\!\!=\!\!&\!\!
-\frac{335267731}{11176704000}\cdot x^{21}+\frac{38742049}{276756480}\cdot x^{19}-\frac{602542}{225225}\cdot x^{17} \\[1.25ex]
\!\!&\!\! \!\!&\!\!
-\frac{169795501}{28828800}\cdot x^{15}+\frac{7838771}{20592}\cdot x^{13}-\frac{15677779}{5280}\cdot x^{11} \\[1.25ex]
\!\!&\!\! \!\!&\!\! +76650\cdot x^{7}-165312\cdot x^{5}-132300\cdot x^{3}+241920\cdot x
\end{array}
\end{array}
$$
and
$$
\begin{array}{l}
\begin{array}{rcl}
P^{-}_{3,3,3,3}(x)
\!\!&\!\!=\!\!&\!\!
\left(-x^{10} -1575\cdot x^{4} -136080\right)\cdot \overline{T}_{\,4\cdot 3+1}^{\,\sin,0}(x) \\[1.25ex]
\!\!&\!\! \!\!&\!\!
\left(-15309\cdot x^{5}-105840\cdot x\right)\cdot \overline{T}_{\,4\cdot 3}^{\,\cos,0}(3 x)
 \\[1.25ex]
\!\!&\!\! \!\!&\!\!
+\left(-\frac{15309}{4}\cdot x^{6}-119070\cdot x^{2}\right)\cdot \overline{T}_{\,4\cdot 3+1}^{\,\sin,0}(3 x)
\\[1.25ex]
\!\!&\!\! \!\!&\!\!
+\left(-\frac{3}{4}\cdot x^{7}-9450\cdot x^{3}\right)\cdot \overline{T}_{\,4\cdot 3}^{\,\cos,0}(x)  \,,
\\[2.25ex]
\end{array} \\[3.00 ex]
\begin{array}{rcl}
P^{-}_{3,3,3,3}(x)
\!\!&\!\!=\!\!&\!\!
-\frac{1}{6227020800}\cdot x^{23}+\frac{1}{39916800}\cdot x^{21}-\frac{1355975527}{1383782400}\cdot x^{19} \\[1.25ex]
\!\!&\!\! \!\!&\!\!
+\frac{4289}{21621600}\cdot x^{17}+\frac{317827231}{28828800}\cdot x^{15}-\frac{8632499}{20592}\cdot x^{13} +\frac{5302897}{1760}\cdot x^{11} \\[1.25ex]
\!\!&\!\! \!\!&\!\!
-76650\cdot x^{7}+165312\cdot x^{5}+132300\cdot x^{3}-241920\cdot x  \,.
\end{array}
\end{array}
$$
Finally, for the MTP function
$$
f(x)=f^{+}(x)+f^{-}(x)
$$
we obtain the concrete {\em downward polynomial approximation}
$$
\begin{array}{rcl}
P(x)
\!\!&\!\!=\!\!&\!\!
P^{+}_{3,3,3,3}(x)+P^{-}_{3,3,3,3}(x)\\[1.5ex]
\!\!&\!\!=\!\!&\!\!
-\frac{1}{6227020800}\cdot x^{23}-\frac{12417313}{413952000}\cdot x^{21}-\frac{581132641}{691891200}\cdot x^{19} \\[1.25ex]
\!\!&\!\! \!\!&\!\!
-\frac{4449211}{1663200}\cdot x^{17}+\frac{21361}{4160}\cdot x^{15}-\frac{424}{11}\cdot x^{13} +\frac{656}{15}\cdot x^{11}
\end{array}
$$
over $\left(0, 1 \right]$, i.e. it holds that
$$
f(x) > P(x)
$$
over $\left(0, 1 \right]$.

\medskip
\noindent
\medskip
\noindent
\textbf{IV} (The final part)
Based on the Sturm theorem, the following inequality
$$
P(x)>0
$$
is true over $\left(0, 1 \right]$.
The stated conclusion for the polynomial function $P$ is correct based on the following facts$:$\\
$1.$
We can conclude, by Sturm theorem, that the polynomial function $P(x)$ has only one zero over the concrete extended segment $[-0.1,1]$
of the initial interval $\left(0,1\right]$.\\
$2.$
Facts
$P\left(-0.1\right) \neq 0$, $P\left(0\right) = 0$ and $P\left(1\right) \neq 0$
are correct. \\
$3.$ The polynomial $P$ is positive at boundary point $1$
{\big (}$P\left(1\right) = 6.77772 \ldots > 0\,${\big )}.
$\!\!$\\[0.1 ex]
Therefore, the following inequality
$$
f(x)>0
$$
is true over  $\left(0,1\right]$. \hfill $\Box$

\bigskip\noindent
{\bf APPENDIX A3} % A3 A3 A3 A3 A3 A3 A3 A3 A3 A3 A3 A3 A3 A3 A3 A3 A3 A3 A3 A3 A3 A3 A3 A3 A3 A3 A3 A3 A3

\medskip\noindent
The initial MTP function is
$$
\begin{array}{rcl}
h(x)
\!\!&\!\!=\!\!&\!\!
\left(-52488 \cdot x^7+816480 \cdot x^5-3810240 \cdot x^3+4354560 \cdot x\right) \cdot \cos^3 x \,                \\[1.25ex]
\!\!&\!\! \!\!&\!\!
 +
\left(-6561\!\cdot\! x^8\!\!+\!244944\!\cdot\! x^6\!\!-\!2041200\!\cdot\! x^4\!\!+\!5080320 \!\cdot\! x^2\!\!-\!1814400\right)\!\cdot\! \sin x \!\cdot\! \cos^2 x\, \\[1.25ex]
\!\!&\!\! \!\!&\!\!
 +
\left(-x^{11} + 39384 \cdot x^7 - 614880 \cdot x^5+2963520 \cdot x^3 - 4354560 \cdot x\right) \cdot \cos x \,      \\[1.25ex]
\!\!&\!\! \!\!&\!\!
 +
\left(1641 \cdot x^8-61488 \cdot x^6+529200 \cdot x^4-1693440 \cdot x^2+1814400\right) \cdot \sin x
\end{array}
$$
and the initial interval is $\SSS = \left[1, \frac{\pi}{2} \right]$.

\break

\medskip
\noindent
After the multiplication by $-1$, we obtain the MTP function
$$
\begin{array}{rcl}
f(x)
\!\!&\!\!=\!\!&\!\!
\left(52488 \cdot x^7 - 816480 \cdot x^5 + 3810240 \cdot x^3-4354560 \cdot x\right) \cdot \cos^3 x \,                 \\[1.25ex]
\!\!&\!\! \!\!&\!\!
 +
%\left(6561x^8 - 244944x^6 + 2041200x^4\ - 5080320x^2 + 1814400\right)\sin x\cos^2 x\, +\\[3.0 ex]
\left(6561 \!\cdot\! x^8\!\!-\!244944 \!\cdot\! x^6\!\!+\!2041200 \!\cdot\! x^4\!\!-\!5080320 \!\cdot\! x^2\!\!+\!1814400\right) \!\cdot\! \sin x \!\cdot\! \cos^2 x\, \\[1.25ex]
\!\!&\!\! \!\!&\!\!
 +
\left(x^{11} - 39384 \cdot x^7 + 614880 \cdot x^5-2963520 \cdot x^3 + 4354560 \cdot x\right) \cdot \cos x \,     \\[1.25ex]
\!\!&\!\! \!\!&\!\!
 +
\left(- 1641 \cdot x^8 + 61488\cdot x^6 - 529200 \cdot x^4 + 1693440 \cdot x^2 - 1814400\right) \cdot \sin x
\,.
\end{array}
$$

\noindent
\underline{\em Automated proof that $f(x) > 0$ for $x \!\in\! \left[ 1, \frac{\pi}{2} \right]$$:$}

\medskip
\noindent
\textbf{I} (Recognition of possible case)
Facts
$f\left(1\right) \neq 0$
and
$f\left(\frac{\pi}{2}\right) \neq 0$
are correct.
The MTP function $f$ is positive at boundary point $1$ {\big (} $\!f\left(1\right) = 27.02986 \ldots  >0$ {\big )} and at boundary point $\frac{\pi}{2}$ {\big (} $\!f\left(\frac{\pi}{2}\right) = 5021.73462 \ldots >0\,${\big )}.
Therefore, it is possible that $f(x) > 0$ over $\left[ 1, \frac{\pi}{2} \right]$.

\medskip
\noindent
\textbf{II} (Transformation of angles)
After the transformation of terms $\cos^m x \cdot \sin^n x$ $(m, n \!\in\! N_0)$
into the sum of sine and cosine functions of multiple angles, in the MTP function $f$,
we obtain
$$
\begin{array}{rcl}
f(x)
\!\!&\!\!=\!\!&\!\!
\left( \frac{6561}{4} \cdot x^8 - 61236 \cdot x^6 + 510300 \cdot x^4 - 1270080 \cdot x^2 + 453600 \right) \cdot \sin 3x   \\[1.25ex]
\!\!&\!\! \!\!&\!\!
+\left( 13122 \cdot x^7 - 204120 \cdot x^5 + 952560 \cdot x^3 - 1088640 \cdot x \right) \cdot \cos 3x  \\[1.25ex]
\!\!&\!\! \!\!&\!\!
+\left( - \frac{3}{4} \cdot x^8 + 252 \cdot x^6 - 18900 \cdot x^4 + 423360 \cdot x^2 - 1360800 \right) \cdot \sin x  \\[1.25ex]
\!\!&\!\! \!\!&\!\!
+\left( x^{11} - 18 \cdot x^7 + 2520 \cdot x^5 - 105840 \cdot x^3 + 1088640 \cdot x \right) \cdot \cos x
 \,.
\end{array}
$$
Then, we consider the previous expression as two separate expressions, the first with positive and the second with negative terms next to sine and cosine functions
$$
\begin{array}{l}
\begin{array}{rcl}
f^{+}(x)
\!\!&\!\!=\!\!&\!\!
\left( \frac{6561}{4} \cdot x^8 + 510300 \cdot x^4 + 453600 \right) \cdot \sin 3x  \\[1.25ex]
\!\!&\!\! \!\!&\!\!
+\left( 13122 \cdot x^7 + 952560 \cdot x^3 \right) \cdot \cos 3x  \\[1.25ex]
\!\!&\!\! \!\!&\!\!
+\left( 252 \cdot x^6 + 423360 \cdot x^2 \right) \cdot \sin x  \\[1.25ex]
\!\!&\!\! \!\!&\!\!
+\left( x^{11} + 2520 \cdot x^5 + 1088640 \cdot x \right) \cdot \cos x
\,, \\[2.25ex]
\end{array} \\[3.00 ex]
\begin{array}{rcl}
f^{-}(x)
\!\!&\!\!=\!\!&\!\!
\left( -61236 \cdot x^6  - 1270080 \cdot x^2  \right) \cdot \sin 3x  \\[1.25ex]
\!\!&\!\! \!\!&\!\!
+\left(- 204120 \cdot x^5 - 1088640 \cdot x \right) \cdot \cos 3x  \\[1.25ex]
\!\!&\!\! \!\!&\!\!
+\left( - \frac{3}{4} \cdot x^8 - 18900 \cdot x^4 - 1360800 \right) \cdot \sin x \\[1.25ex]
\!\!&\!\! \!\!&\!\!
+\left( -18 \cdot x^7 - 105840 \cdot x^3 \right) \cdot \cos x \,.
\end{array}
\end{array}
$$

\medskip\noindent
\textbf{III} (Determination of downward rational polynomial approximation)
After substitution of sine and cosine functions by appropriate (downward or upward)
polynomial approximations, we obtain downward polynomial approximations of $f^{+}(x)$ and $f^{-}(x)$
respectively
$$
\begin{array}{l}
\begin{array}{rcl}
P^{+}_{i_{0},i_{1},i_{2},i_{3}}(x)
\!\!&\!\!=\!\!&\!\!
\left( \frac{6561}{4} \cdot x^8 + 510300 \cdot x^4 + 453600 \right) \cdot \underline{T}_{4\cdot i_{0}+3}^{\sin,0}(3 x)  \\[1.25ex]
\!\!&\!\! \!\!&\!\!
+\left( 13122 \cdot x^7 + 952560 \cdot x^3 \right) \cdot \underline{T}_{4\cdot i_{1}+2}^{\cos,0}(3 x)  \\[1.25ex]
\!\!&\!\! \!\!&\!\!
+\left( 252 \cdot x^6 + 423360 \cdot x^2 \right) \cdot \underline{T}_{4\cdot i_{2}+3}^{\sin,0}(x)  \\[1.25ex]
\!\!&\!\! \!\!&\!\!
+\left( x^{11} + 2520 \cdot x^5 + 1088640 \cdot x \right) \cdot \underline{T}_{4\cdot i_{3}+2}^{\cos,0}(x) \,,
\\[2.25ex]
\end{array} \\[3.00 ex]
\begin{array}{rcl}
P^{-}_{i_{4},i_{5},i_{6},i_{7}}(x)
\!\!&\!\!=\!\!&\!\!
\left( -61236 \cdot x^6  - 1270080 \cdot x^2  \right) \cdot \overline{T}_{4\cdot i_{4}+1}^{\sin,0}(3 x)   \\[1.25ex]
\!\!&\!\! \!\!&\!\!
+\left( -204120 \cdot x^5 - 1088640 \cdot x \right) \cdot \overline{T}_{4\cdot i_{5}}^{\cos,0}(3 x)  \\[1.25ex]
\!\!&\!\! \!\!&\!\!
+\left(  -\frac{3}{4} \cdot x^8 - 18900 \cdot x^4 - 1360800 \right) \cdot \overline{T}_{4\cdot i_{6}+1}^{\sin,0}(x)  \\[1.25ex]
\!\!&\!\! \!\!&\!\!
+\left( -18 \cdot x^7- 105840 \cdot x^3 \right) \cdot \overline{T}_{4\cdot i_{7}}^{\cos,0}(x)
\,.
\end{array}
\end{array}
$$
For concrete indices $\left(i_{0},i_{1},i_{2},i_{3},i_{4},i_{5},i_{6},i_{7}\right)=\left(3,4,1,2,4,4,2,2\right)$, we obtain
$$
\begin{array}{l}
\begin{array}{rcl}
P^{+}_{3,4,1,2}(x)
\!\!&\!\!=\!\!&\!\!
\left( \frac{6561}{4} \cdot x^8 + 510300 \cdot x^4 + 453600 \right) \cdot \underline{T}_{4\cdot 3+3}^{\sin,0}(3 x)  \\[1.25ex]
\!\!&\!\! \!\!&\!\!
+\left( 13122 \cdot x^7 + 952560 \cdot x^3 \right) \cdot \underline{T}_{4\cdot 4+2}^{\cos,0}(3 x)  \\[1.25ex]
\!\!&\!\! \!\!&\!\!
+\left( 252 \cdot x^6 + 423360 \cdot x^2 \right) \cdot \underline{T}_{4\cdot 1+3}^{\sin,0}(x) \\[1.25ex]
\!\!&\!\! \!\!&\!\!
+\left( x^{11} + 2520 \cdot x^5 + 1088640 \cdot x \right) \cdot \underline{T}_{4\cdot 2+2}^{\cos,0}(x) \,,
\\[2.25ex]
\end{array} \\[3.00 ex]
\begin{array}{rcl}
P^{+}_{3,4,1,2}(x)
\!\!&\!\!=\!\!&\!\!
- \frac{387420489}{487911424000} \cdot x^{25}+ \frac{129140163}{14350336000} \cdot x^{23} - \frac{2208298489}{6175128960} \cdot x^{21}
\\[1.25ex]
\!\!&\!\! \!\!&\!\!
+ \frac{129141043}{35481600} \cdot x^{19} - \frac{760514951}{16473600} \cdot x^{17} + \frac{607555181}{2882880} \cdot x^{15}
\\[1.25ex]
\!\!&\!\! \!\!&\!\!
+ \frac{53310093}{22880} \cdot x^{13}- \frac{13876823}{440} \cdot x^{11} + 39372 \cdot x^9+735840 \cdot x^{7}
\\[1.25ex]
\!\!&\!\! \!\!&\!\!
-1859760 \cdot x^{5} -1209600 \cdot x^{3} + 2449440 \cdot x
\end{array}
\end{array}
$$
and
$$
\begin{array}{l}
\begin{array}{rcl}
P^{-}_{4,4,2,2}(x)
\!\!&\!\!=\!\!&\!\!
\left( -61236 \cdot x^6 - 1270080 \cdot x^2  \right) \cdot \overline{T}_{4\cdot 4+1}^{\sin,0}(3 x)   \\[1.25ex]
\!\!&\!\! \!\!&\!\!
+\left( -204120 \cdot x^5 - 1088640 \cdot x \right) \cdot \overline{T}_{4\cdot 4}^{\cos,0}(3 x) \\[1.25ex]
\!\!&\!\! \!\!&\!\!
+\left(  -\frac{3}{4} \cdot x^8 - 18900 \cdot x^4 - 1360800 \right) \cdot \overline{T}_{4\cdot 2+1}^{\sin,0}(x) \\[1.25ex]
\!\!&\!\! \!\!&\!\!
+\left( -18 \cdot x^7- 105840 \cdot x^3 \right) \cdot \overline{T}_{4\cdot 2}^{\cos,0}(x)
\,, \\[2.25ex]
\end{array} \\[3.00 ex]
\begin{array}{rcl}
P^{-}_{4,4,2,2}(x)
\!\!&\!\!=\!\!&\!\!
- \frac{387420489}{17425408000} \cdot x^{23} + \frac{129140163}{512512000} \cdot x^{21} - \frac{43046721}{8712704} \cdot x^{19}
\\[1.25ex]
\!\!&\!\! \!\!&\!\!
 + \frac{19715397503}{345945600} \cdot x^{17} - \frac{127545983}{480480} \cdot x^{15} - \frac{531449}{240} \cdot x^{13}
\\[1.25ex]
\!\!&\!\! \!\!&\!\!
 +\frac{314933}{10} \cdot x^{11} - 39372 \cdot x^9  - 735840 \cdot x^7
\\[1.25ex]
\!\!&\!\! \!\!&\!\!
+ 1859760 \cdot x^5 + 1209600 \cdot x^3 - 2449440 \cdot x
 \,.
\end{array}
\end{array}
$$
Finally, for the MTP function\\
$$
f(x)
=
f^{+}(x)+f^{-}(x)
$$
we obtain the concrete {\em downward polynomial approximation}
$$
\begin{array}{rcl}
P(x)
\!\!&\!\!=\!\!&\!\!
P^{+}_{3,4,1,2}(x)+P^{-}_{4,4,2,2}(x) \\[1.5 ex]
\!\!&\!\!=\!\!&\!\!
-\frac{387420489}{487911424000} \cdot x^{25} - \frac{129140163}{9758228480} \cdot x^{23} - \frac{521856760043}{4940103168000} \cdot x^{21}
\\[1.25ex]
\!\!&\!\! \!\!&\!\!
-\frac{10201878397}{7841433600} \cdot x^{19} + \frac{936145883}{86486400} \cdot x^{17} - \frac{2048321}{37440} \cdot x^{15}
\\[1.25ex]
\!\!&\!\! \!\!&\!\!
+\frac{1587173}{13728} \cdot x^{13} - \frac{19771}{440} \cdot x^{11}
\end{array}
$$
over $\left[1,\frac{\pi}{2}\right]$, i.e. it holds that
$$
f(x) > P(x)
$$
over $\left[1,\frac{\pi}{2}\right]$.

\medskip
\noindent
\textbf{IV} (The final part)
Based on the Sturm theorem, the following inequality
$$
P(x)>0
$$
is true over $\left[1,\frac{\pi}{2}\right]$.
The stated conclusion for the polynomial function $P$ is correct based on the following facts$:$\\
$1.$
We can conclude, by Sturm theorem, that the polynomial function $P(x)$ does not have zero over the concrete extended segment $[1, 1.58]$
of the initial interval $\left[1,\frac{\pi}{2}\right]$.\\
$2.$
Facts
$P\left(1\right) \neq 0$ and $P\left(1.58 \right) \neq 0$
are correct.\\
$3.$ The polynomial $P$ is positive at boundary point  $\frac{\pi}{2}$
{\big (}$P\left(\frac{\pi}{2}\right) = 1228.02881 \ldots > 0\,${\big )}.\\[0.1 ex]
$\!\!$Therefore, the following inequality
$$
f(x)>0
$$
is true over  $\left[1,\frac{\pi}{2}\right]$. \hfill $\Box$

\bigskip\noindent
{\bf APPENDIX A4} % A4 A4 A4 A4 A4 A4 A4 A4 A4 A4 A4 A4 A4 A4 A4 A4 A4 A4 A4 A4 A4 A4 A4 A4 A4 A4 A4 A4 A4

\medskip\noindent
The initial MTP function is
$$
\begin{array}{rcl}
f(x)=x^{5} -360\cdot x+\left(-30\cdot x^{2} +630\right)\cdot\sin x -270\cdot x\cdot\cos x
\end{array}
$$
and the initial interval is $\SSS = \left(0, \frac{\pi}{2} \right)$.

\medskip\noindent
\underline{\em Automated proof that $f(x) > 0$ for $x \!\in\! \left(0, \frac{\pi}{2} \right)$$:$}

\medskip
\noindent
\textbf{I} (Recognition of possible case)
Facts
$f\left(0\right) = 0$
and
$f\left(\frac{\pi}{2}\right) \neq 0$
are correct.
The MTP function $f$ is positive at boundary point $\frac{\pi}{2}$ {\big (} $\!f\left(\frac{\pi}{2}\right) = 0.054404 \ldots >0${\big )}.
Therefore, it is possible that $f(x) > 0$ over $\left(0,\frac{\pi}{2}\right)$.

\medskip
\noindent
\textbf{II} (Transformation of angles)
After the transformation of terms $\cos^m x \cdot \sin^n x$ $(m, n \!\in\! N_0)$
into the sum of sine and cosine functions of multiple angles, in the MTP function $f$,
we obtain
$$
\begin{array}{rcl}
f(x)=\left(-30\cdot x^{2} +630\right)\cdot\sin x -270\cdot x\cdot\cos x +x^{5} -360\cdot x \,.
\end{array}
$$
Then, we consider the previous expression as two separate expressions, the first with positive and the second with negative terms next to sine and cosine functions
$$
\begin{array}{l}
\begin{array}{rcl}
f^{+}(x)
\!\!&\!\!=\!\!&\!\!
630\cdot\sin x +x^{5} \,,
\end{array} \\[2.50 ex]
\begin{array}{rcl}
f^{-}(x)
\!\!&\!\!=\!\!&\!\!
-360\cdot x-270\cdot x\cdot\cos x -30\cdot x^{2}\cdot\sin x \,.
\end{array}
\end{array}
$$

\medskip\noindent
\textbf{III} (Determination of downward rational polynomial approximation)
After substitution of sine and cosine functions by appropriate (downward or upward)
polynomial approximations, we obtain downward polynomial approximations of $f^{+}(x)$ and $f^{-}(x)$
respectively
$$
\begin{array}{l}
\begin{array}{rcl}
P^{+}_{i_{0}}(x)=630\cdot \underline{T}_{4\cdot i_{0}+3}^{\sin,0}(x)+x^{5} \,,
\end{array} \\[2.50 ex]
\begin{array}{rcl}
P^{-}_{i_{1},i_{2}}(x)=-360\cdot x-270\cdot x\cdot \overline{T}_{4\cdot i_{1}}^{\cos,0}(x)-30\cdot x^{2}\cdot \overline{T}_{4\cdot i_{2}+1}^{\sin,0}(x) \,.
\end{array}
\end{array}
$$
For concrete indices $\left(i_{0},i_{1},i_{2}\right)=\left(2,2,2\right)$, we obtain
$$
\begin{array}{l}
\begin{array}{rcl}
P^{+}_{2}(x)=630\cdot \underline{T}_{4\cdot 2+3}^{\sin,0}(x)+x^{5}  \,,
\end{array} \\[2.50 ex]
\begin{array}{rcl}
P^{+}_{2}(x)=-\frac{1}{63360}\cdot x^{11}+\frac{1}{576}\cdot x^{9}-\frac{1}{8}\cdot x^{7}+\frac{25}{4}\cdot x^{5}-105\cdot x^{3} + 630\cdot x
\end{array}
\end{array}
$$
and
$$
\begin{array}{l}
\begin{array}{rcl}
P^{-}_{2,2}(x)=
-360\cdot x-270\cdot x\cdot \overline{T}_{4\cdot 2}^{\cos,0}(x)-30\cdot x^{2}\cdot \overline{T}_{4\cdot 2+1}^{\sin,0}(x)  \,,
\end{array} \\[2.50 ex]
\begin{array}{rcl}
P^{-}_{2,2}(x)=-\frac{1}{12096}\cdot x^{11}-\frac{1}{1344}\cdot x^{9}+\frac{1}{8}\cdot x^{7}-\frac{25}{4}\cdot x^{5}+105\cdot x^{3} -  630\cdot x \,.
\end{array}
\end{array}
$$
Finally, for the MTP function
$$
f(x)
=
f^{+}(x)+f^{-}(x)
$$
we obtain the concrete {\em downward polynomial approximation}
$$
\begin{array}{rcl}
P(x)
\!\!&\!\!=\!\!&\!\!
P^{+}_{2}(x)+P^{-}_{2,2}(x) \\[1.5 ex]
\!\!&\!\!=\!\!&\!\!
-\frac{131}{1330560}\cdot x^{11}+\frac{1}{1008}\cdot x^{9}
\end{array}
$$
over $\left(0,\frac{\pi}{2}\right)$, i.e. it holds that
$$
f(x) > P(x)
$$
over $\left(0,\frac{\pi}{2}\right)$.

\medskip
\noindent
\textbf{IV} (The final part)
Based on the Sturm theorem, the following inequality
$$
P(x)>0
$$
is true over $\left(0,\frac{\pi}{2}\right)$.
The stated conclusion for the polynomial function $P$ is correct based on the following facts$:$\\
$1.$
We can conclude, by Sturm theorem, that the polynomial function $P(x)$ has only one zero over the concrete extended segment $[-0.1, 1.58]$
of the initial interval $\left(0,\frac{\pi}{2}\right)$.\\
$2.$
Facts
$P\left(-0.1\right) \neq 0$, $P\left(0\right) = 0$ and $P\left(1.58 \right) \neq 0$
are correct.\\
$3.$ The polynomial $P$ is positive at boundary point $\frac{\pi}{2}$
{\big (}$P\left(\frac{\pi}{2}\right) = 0.043615 \ldots > 0\,${\big )}.
$\!\!$\\[0.1 ex]
Therefore, the following inequality
$$
f(x)>0
$$
is true over $\left(0,\frac{\pi}{2}\right)$. \hfill $\Box$

\bigskip\noindent
{\bf APPENDIX A5} % A5 A5 A5 A5 A5 A5 A5 A5 A5 A5 A5 A5 A5 A5 A5 A5 A5 A5 A5 A5 A5 A5 A5 A5 A5 A5 A5 A5 A5

\medskip\noindent
The initial MTP function is
$$
\begin{array}{rcl}
f(x)
\!\!&\!\!=\!\!&\!\!
\left(x^{8} -21\cdot x^{6} +630\cdot x^{4} -7560\cdot x^{2} +15120\right)\cdot\sin x                 \\[1.25ex]
\!\!&\!\! \!\!&\!\!
 +
\left(3\cdot x^{7} -126\cdot x^{5} +2520\cdot x^{3} -15120\cdot x\right)\cdot\cos x
\end{array}
$$
and the initial interval is $\SSS = \left(0, \frac{\pi}{2} \right)$.

\break

% \medskip
\noindent
\underline{\em Automated proof that $f(x) > 0$ for $x \!\in\! \left(0, \frac{\pi}{2} \right)$$:$}

\medskip
\noindent
\textbf{I} (Recognition of possible case)
Facts
$f\left(0\right) = 0$
and
$f\left(\frac{\pi}{2}\right) \neq 0$
are correct.
The MTP function $f$ is positive at boundary point $\frac{\pi}{2}$ {\big (} $\!f\left(\frac{\pi}{2}\right) =  23.53938 \ldots >0\,${\big )}.
Therefore, it is possible that $f(x) > 0$ over $\left(0,\frac{\pi}{2}\right)$.

\medskip
\noindent
\textbf{II} (Transformation of angles)
After the transformation of terms $\cos^m x \cdot \sin^n x$ $(m, n \!\in\! N_0)$
into the sum of sine and cosine functions of multiple angles, in the MTP function $f$,
we obtain
$$
\begin{array}{rcl}
f(x)
\!\!&\!\!=\!\!&\!\!
\left(x^{8} -21\cdot x^{6} +630\cdot x^{4} -7560\cdot x^{2} +15120\right)\cdot\sin x  \\[1.25ex]
\!\!&\!\! \!\!&\!\!
+ \left(3\cdot x^{7} -126\cdot x^{5} +2520\cdot x^{3} -15120\cdot x\right)\cdot\cos x \,.
 \end{array}
$$
Then, we consider the previous expression as two separate expressions, the first with positive and the second with negative terms next to sine and cosine functions
$$
\begin{array}{l}
\begin{array}{rcl}
f^{+}(x)
\!\!&\!\!=\!\!&\!\!
\left(x^{8} +630\cdot x^{4} +15120\right)\cdot\sin x +\left(3\cdot x^{7} +2520\cdot x^{3}\right)\cdot\cos x \,,
\end{array} \\[2.50 ex]
\begin{array}{rcl}
f^{-}(x)
\!\!&\!\!=\!\!&\!\!
\left(-126\cdot x^{5} -15120\cdot x\right)\cdot\cos x +\left(-21\cdot x^{6} -7560\cdot x^{2}\right)\cdot\sin x \,.
\end{array}
\end{array}
$$

\medskip\noindent
\textbf{III} (Determination of downward rational polynomial approximation)
After substitution of sine and cosine functions by appropriate (downward or upward)
polynomial approximations, we obtain downward polynomial approximations of $f^{+}(x)$ and $f^{-}(x)$
respectively
$$
\begin{array}{l}
\begin{array}{rcl}
P^{+}_{i_{0},i_{1}}(x)=\left(x^{8} +630\cdot x^{4} +15120\right)\cdot \underline{T}_{4\cdot i_{0}+3}^{\sin,0}(x)+\left(3\cdot x^{7} +2520\cdot x^{3}\right)\cdot \underline{T}_{4\cdot i_{1}+2}^{\cos,0}(x) \,,
\end{array} \\[2.50 ex]
\begin{array}{rcl}
P^{-}_{i_{2},i_{3}}(x)=\left(-126\cdot x^{5} -15120\cdot x\right)\cdot \overline{T}_{4\cdot i_{2}}^{\cos,0}(x)+\left(-21\cdot x^{6} -7560\cdot x^{2}\right)\cdot \overline{T}_{4\cdot i_{3}+1}^{\sin,0}(x) \,.
\end{array}
\end{array}
$$
For concrete indices $\left(i_{0},i_{1},i_{2},i_{3}\right)=\left(1,1,2,2\right)$, we obtain
$$
\begin{array}{l}
\begin{array}{rcl}
P^{+}_{1,1}(x)=\left(x^{8} +630\cdot x^{4} +15120\right)\cdot \underline{T}_{4\cdot 1+3}^{\sin,0}(x)+\left(3\cdot x^{7} +2520\cdot x^{3}\right)\cdot \underline{T}_{4\cdot 1+2}^{\cos,0}(x)  \,,
\end{array} \\[2.50 ex]
\begin{array}{rcl}
P^{+}_{1,1}(x)=-\frac{1}{5040}\cdot x^{15}+\frac{1}{240}\cdot x^{13}-\frac{1}{6}\cdot x^{11}+\frac{5}{4}\cdot x^{9}-504\cdot x^{5} +15120\cdot x
\end{array}
\end{array}
$$
and
$$
\begin{array}{l}
\begin{array}{rcl}
P^{-}_{2,2}(x)=
\left(-126\cdot x^{5} -15120\cdot x\right)\cdot \overline{T}_{4\cdot 2}^{\cos,0}(x)+\left(-21\cdot x^{6} -7560\cdot x^{2}\right)\cdot \overline{T}_{4\cdot 2+1}^{\sin,0}(x)  \,,
\end{array} \\[2.50 ex]
\begin{array}{rcl}
P^{-}_{2,2}(x)=-\frac{1}{17280}\cdot x^{15}+\frac{1}{960}\cdot x^{13}-\frac{1}{48}\cdot x^{11}-\frac{5}{8}\cdot x^{9}+504\cdot x^{5} -15120\cdot x \,.
\end{array}
\end{array}
$$
Finally, for the MTP function
$$
f(x)
=
f^{+}(x)+f^{-}(x)
$$
we obtain the concrete {\em downward polynomial approximation}
$$
\begin{array}{rcl}
P(x)
\!\!&\!\!=\!\!&\!\!
P^{+}_{1,1}(x)+P^{-}_{2,2}(x) \\[1.5 ex]
\!\!&\!\!=\!\!&\!\!
-\frac{31}{120960}\cdot x^{15}+\frac{1}{192}\cdot x^{13}-\frac{3}{16}\cdot x^{11}+\frac{5}{8}\cdot x^{9}
\end{array}
$$
over $\left(0,\frac{\pi}{2}\right)$, i.e. it holds that
$$
f(x) > P(x)
$$
over $\left(0,\frac{\pi}{2}\right)$.

\medskip
\noindent
\textbf{IV} (The final part)
Based on the Sturm theorem, the following inequality
$$
P(x)>0
$$
is true over $\left(0,\frac{\pi}{2}\right)$.
The stated conclusion for the polynomial function $P$ is correct based on the following facts$:$\\
$1.$
We can conclude, by Sturm theorem, that the polynomial function $P(x)$ has only one zero over the concrete extended segment $[-0.1, 1.58]$
of the initial interval $\left(0,\frac{\pi}{2}\right)$.\\
$2.$
Facts
$P\left(-0.1\right) \neq 0$, $P\left(0\right) = 0$ and $P\left(1.58 \right) \neq 0$
are correct.\\
$3.$ The polynomial $P$ is positive at boundary point $\frac{\pi}{2}$
{\big (}$P\left(\frac{\pi}{2}\right) =11.074847 \ldots > 0\,${\big )}.
$\!\!$\\[0.1 ex]
Therefore, the following inequality
$$
f(x)>0
$$
is true over  $\left(0,\pi/2\right)$. \hfill $\Box$

\bigskip\noindent
{\bf APPENDIX A6} % A6 A6 A6 A6 A6 A6 A6 A6 A6 A6 A6 A6 A6 A6 A6 A6 A6 A6 A6 A6 A6 A6 A6 A6 A6 A6 A6 A6 A6

\medskip\noindent
The initial MTP function is
$$
\begin{array}{rcl}
f(x)
\!\!&\!\!=\!\!&\!\!
3\cdot x^{5} -20\cdot x-140\cdot\sin x \cdot\cos x +\left(30\cdot x^{2} -70\right)\cdot\sin x         \\[1.25ex]
\!\!&\!\! \!\!&\!\!
 +
40\cdot x\cdot\cos^{2} x +190\cdot x\cdot\cos x
\end{array}
$$
and the initial interval is $\SSS = \left(0, \frac{\pi}{2} \right)$.

\medskip\noindent
\underline{\em Automated proof that $f(x) > 0$ for $x \!\in\! \left(0, \frac{\pi}{2} \right)$$:$}

\medskip
\noindent
\textbf{I} (Recognition of possible case)
Facts
$f\left(0\right) = 0$
and
$f\left(\frac{\pi}{2}\right) \neq 0$
are correct.
The MTP function $f$ is positive at boundary point $\frac{\pi}{2}$ {\big (} $\!f\left(\frac{\pi}{2}\right) = 1.29545 \ldots >0${\big )}.
Therefore, it is possible that $f(x) > 0$ over $\left(0,\frac{\pi}{2}\right)$.

\medskip
\noindent
\textbf{II} (Transformation of angles)
After the transformation of terms $\cos^m x \cdot \sin^n x$ $(m, n \!\in\! N_0)$
into the sum of sine and cosine functions of multiple angles, in the MTP function $f$,
we obtain
$$
f(x)
=
-70\cdot\sin 2 x +20\cdot x\cdot\cos 2 x +\left(30\cdot x^{2} -70\right)\cdot\sin x +190\cdot x\cdot\cos x+3\cdot x^{5}
\,.
$$
Then, we consider the previous expression as two separate expressions, the first with positive and the second with negative terms next to sine and cosine functions
$$
\begin{array}{l}
\begin{array}{rcl}
f^{+}(x)
\!\!&\!\!=\!\!&\!\!
190\cdot x\cdot\cos x +20\cdot x\cdot\cos 2 x +30\cdot x^{2}\cdot\sin x +3\cdot x^{5}\,,
\end{array} \\[2.50 ex]
\begin{array}{rcl}
f^{-}(x)
\!\!&\!\!=\!\!&\!\!
-70\cdot\sin x -70\cdot\sin 2 x \,.
\end{array}
\end{array}
$$

\medskip\noindent
\textbf{III} (Determination of downward rational polynomial approximation)
After substitution of sine and cosine functions by appropriate (downward or upward)
polynomial approximations, we obtain downward polynomial approximations of $f^{+}(x)$ and $f^{-}(x)$
respectively
$$
\begin{array}{l}
\begin{array}{rcl}
\!\!\!\!\!\!P^{+}_{i_{0},i_{1},i_{2}}(x)
\!\!&\!\!=\!\!&\!\!
190\cdot x\cdot \underline{T}_{4\cdot i_{0}+2}^{\cos,0}(x)+20\cdot x\cdot \underline{T}_{4\cdot i_{1}+2}^{\cos,0}(2 x) +30\cdot x^{2}\cdot \underline{T}_{4\cdot i_{2}+3}^{\sin,0}(x)+3\cdot x^{5}
\,,
\end{array} \\[2.50 ex]
\begin{array}{rcl}
\!\!\!\!\!\!P^{-}_{i_{3},i_{4}}(x)
\!\!&\!\!=\!\!&\!\!
-70\cdot \overline{T}_{4\cdot i_{3}+1}^{\sin,0}(x)-70\cdot \overline{T}_{4\cdot i_{4}+1}^{\sin,0}(2 x)
\,.
\end{array}
\end{array}
$$
For concrete indices $\left(i_{0},i_{1},i_{2},i_{3},i_{4}\right)=\left(1, 2, 1, 2, 2\right)$, we obtain
$$
\begin{array}{l}
\begin{array}{rcl}
P^{+}_{1,2,1}(x)
\!\!&\!\!=\!\!&\!\!
190\cdot x\cdot \underline{T}_{4\cdot 1+2}^{\cos,0}(x)+20\cdot x\cdot \underline{T}_{4\cdot 2+2}^{\cos,0}(2 x)+30\cdot x^{2}\cdot \underline{T}_{4\cdot 1+3}^{\sin,0}(x)+3\cdot x^{5}
\,,
\end{array} \\[2.50 ex]
\begin{array}{rcl}
P^{+}_{1,2,1}(x)
\!\!&\!\!=\!\!&\!\!
-\frac{16}{2835}\cdot x^{11}+\frac{61}{504}\cdot x^{9}-\frac{43}{24}\cdot x^{7}+\frac{77}{4}\cdot x^{5}-105\cdot x^{3}+210\cdot x
\end{array}
\end{array}
$$
and
$$
\begin{array}{l}
\begin{array}{rcl}
P^{-}_{2,2}(x)
\!\!&\!\!=\!\!&\!\!
-70\cdot \overline{T}_{4\cdot 2+1}^{\sin,0}(x)-70\cdot \overline{T}_{4\cdot 2+1}^{\sin,0}(2 x)\,,
\end{array} \\[2.50 ex]
\begin{array}{rcl}
P^{-}_{2,2}(x)
\!\!&\!\!=\!\!&\!\!
-\frac{19}{192}\cdot x^{9}+\frac{43}{24}\cdot x^{7}-\frac{77}{4}\cdot x^{5}+105\cdot x^{3}-210\cdot x\,.
\end{array}
\end{array}
$$
Finally, for the MTP function
$$
f(x)
=
f^{+}(x)+f^{-}(x)
$$
we obtain the concrete {\em downward polynomial approximation}
$$
\begin{array}{rcl}
P(x)
\!\!&\!\!=\!\!&\!\!
P^{+}_{1,2,1}(x)+P^{-}_{2,2}(x) \\[1.5 ex]
\!\!&\!\!=\!\!&\!\!
-\frac{16}{2835}\cdot x^{11}+\frac{89}{4032}\cdot x^{9}
\end{array}
$$
over $\left(0,\frac{\pi}{2}\right)$, i.e. it holds that
$$
f(x) > P(x)
$$
over $\left(0,\frac{\pi}{2}\right)$.

\medskip
\noindent
\textbf{IV} (The final part)
Based on the Sturm theorem, the following inequality
$$
P(x)>0
$$
is true over $\left(0,\frac{\pi}{2}\right)$.
The stated conclusion for the polynomial function $P$ is correct based on the following facts$:$\\
$1.$
We can conclude, by Sturm theorem, that the polynomial function $P(x)$ has only one zero over the concrete extended segment $[-0.1, 1.58]$
of the initial interval $\left(0,\frac{\pi}{2}\right)$.\\
$2.$
Facts
$P\left(-0.1\right) \neq 0$, $P\left(0\right) = 0$ and $P\left(1.58 \right) \neq 0$
are correct.\\
$3.$ The polynomial $P$ is positive at boundary point $\frac{\pi}{2}$
{\big (}$P\left(\frac{\pi}{2}\right) = 0.47438 \ldots > 0\,${\big )}.
$\!\!$\\[0.1 ex]
Therefore, the following inequality
$$
f(x)>0
$$
is true over $\left(0,\frac{\pi}{2}\right)$. \hfill $\Box$

\bigskip\noindent
{\bf APPENDIX A7} % A7 A7 A7 A7 A7 A7 A7 A7 A7 A7 A7 A7 A7 A7 A7 A7 A7 A7 A7 A7 A7 A7 A7 A7 A7 A7 A7 A7 A7

\medskip\noindent
The initial MTP function is
$$
\begin{array}{rcl}
f(x)=\!\!&\!\!=\!\!&\!\!
16\cdot x^{6}+135\cdot x\cdot\sin x \cdot\cos^{2} x+\left(90\cdot x^{3} -360\cdot x\right)\cdot\sin x \\[1.25ex]
\!\!&\!\! \!\!&\!\!
 +360\cdot\cos^{3} x+\left(585\cdot x^{2} -360\right)\cdot\cos x
\end{array}
$$
and the initial interval is $\SSS = \left(0, \frac{\pi}{2} \right)$.

\medskip\noindent
\underline{\em Automated proof that $f(x) > 0$ for $x \!\in\! \left(0, \frac{\pi}{2} \right)$$:$}

\medskip
\noindent
\textbf{I} (Recognition of possible case)
Facts
$f\left(0\right) = 0$
and
$f\left(\frac{\pi}{2}\right) \neq 0$
are correct.
The MTP function $f$ is positive at boundary point $\frac{\pi}{2}$ {\big (} $\!f\left(\frac{\pi}{2}\right) = 23.68123 \ldots >0\,${\big )}.
Therefore, it is possible that $f(x) > 0$ over $\left(0,\frac{\pi}{2}\right)$.

\medskip
\noindent
\textbf{II} (Transformation of angles)
After the transformation of terms $\cos^m x \cdot \sin^n x$ $(m, n \!\in\! N_0)$
into the sum of sine and cosine functions of multiple angles, in the MTP function $f$,
we obtain
$$
\begin{array}{rcl}
f(x)=\!\!&\!\!=\!\!&\!\!
90\cdot\cos 3x +\left(90\cdot x^{3} -\frac{1305}{4}\cdot x\right)\cdot\sin x +\left(585\cdot x^{2} -90\right)\cdot\cos x  \\[1.25ex]
\!\!&\!\! \!\!&\!\!
+16\cdot x^{6}+\frac{135}{4}\cdot x\cdot\sin 3 x
\,.
 \end{array}
$$
Then, we consider the previous expression as two separate expressions, the first with positive and the second with negative terms next to sine and cosine functions
$$
\begin{array}{l}
\begin{array}{rcl}
f^{+}(x)
\!\!&\!\!=\!\!&\!\!
90\cdot\cos 3 x +\frac{135}{4}\cdot x\cdot\sin 3 x +585\cdot x^{2}\cdot\cos x +90\cdot x^{3}\cdot\sin x+16\cdot x^{6} \,,
\end{array} \\[2.50 ex]
\begin{array}{rcl}
f^{-}(x)
\!\!&\!\!=\!\!&\!\!
-90\cdot\cos x -\frac{1305}{4}\cdot x\cdot\sin x  \,.
\end{array}
\end{array}
$$

\medskip\noindent
\textbf{III} (Determination of downward rational polynomial approximation)
After substitution of sine and cosine functions by appropriate (downward or upward)
polynomial approximations, we obtain downward polynomial approximations of $f^{+}(x)$ and $f^{-}(x)$
respectively
$$
\begin{array}{l}
\begin{array}{rcl}
P^{+}_{i_{0},i_{1},i_{2},i_{3}}(x)
\!\!&\!\!=\!\!&\!\!
90\cdot \underline{T}_{4\cdot i_{0}+2}^{\cos,0}(3 x)+\frac{135}{4}\cdot x\cdot \underline{T}_{4\cdot i_{1}+3}^{\sin,0}(3 x)   \\[1.25ex]
\!\!&\!\! \!\!&\!\!
+585\cdot x^{2}\cdot \underline{T}_{4\cdot i_{2}+2}^{\cos,0}(x)+90\cdot x^{3}\cdot \underline{T}_{4\cdot i_{3}+3}^{\sin,0}(x)+16\cdot x^{6}
\,,
\\[2.25ex]
\end{array} \\[3.00 ex]
\begin{array}{rcl}
P^{-}_{i_{4},i_{5}}(x)
\!\!&\!\!=\!\!&\!\!
-90\cdot \overline{T}_{4\cdot i_{4}}^{\cos,0}(x)-\frac{1305}{4}\cdot x\cdot \overline{T}_{4\cdot i_{5}+1}^{\sin,0}(x)
\,.
\end{array}
\end{array}
$$
For concrete indices $\left(i_{0},i_{1},i_{2},i_{3},i_{4},i_{5}\right)=\left(2, 3, 1, 1, 2, 2\right)$, we obtain
$$
\begin{array}{l}
\begin{array}{rcl}
P^{+}_{2,3,1,1}(x)
\!\!&\!\!=\!\!&\!\!
90\cdot \underline{T}_{4\cdot 2+2}^{\cos,0}(3 x)+\frac{135}{4}\cdot x\cdot \underline{T}_{4\cdot 3+3}^{\sin,0}(3 x)   \\[1.25ex]
\!\!&\!\! \!\!&\!\!
+
585\cdot x^{2}\cdot \underline{T}_{4\cdot 1+2}^{\cos,0}(x)+90\cdot x^{3}\cdot \underline{T}_{4\cdot 1+3}^{\sin,0}(x)+16\cdot x^{6}
\,,
\\[2.25ex]
\end{array} \\[3.00 ex]
\begin{array}{rcl}
P^{+}_{2,3,1,1}(x)
\!\!&\!\!=\!\!&\!\!
-\frac{531441}{1435033600}\cdot x^{16}+\frac{177147}{20500480}\cdot x^{14}-\frac{59049}{394240}\cdot x^{12} + \frac{6241}{17920}\cdot x^{10}   \\[1.25ex]
\!\!&\!\! \!\!&\!\!
-\frac{1}{16}\cdot x^{8}+\frac{83}{32}\cdot x^{6}-\frac{405}{8}\cdot x^{4}+\frac{1125}{4}\cdot x^{2} + 90
\end{array}
\end{array}
$$
and
$$
\begin{array}{l}
\begin{array}{rcl}
P^{-}_{2,2}(x)
\!\!&\!\!=\!\!&\!\!
-90\cdot \overline{T}_{4\cdot 2}^{\cos,0}(x)-\frac{1305}{4}\cdot x\cdot \overline{T}_{4\cdot 2+1}^{\sin,0}(x) \,,
\end{array} \\[2.50 ex]
\begin{array}{rcl}
P^{-}_{2,2}(x)
\!\!&\!\!=\!\!&\!\!
-\frac{29}{32256}\cdot x^{10}+\frac{1}{16}\cdot x^{8}-\frac{83}{32}\cdot x^{6}+\frac{405}{8}\cdot x^{4}-\frac{1125}{4}\cdot x^{2} - 90\,.
\end{array}
\end{array}
$$
Finally, for the MTP function
$$
f(x)
=
f^{+}(x)+f^{-}(x)
$$
we obtain the concrete {\em downward polynomial approximation}
$$
\begin{array}{rcl}
P(x)
\!\!&\!\!=\!\!&\!\!
P^{+}_{2,3,1,1}(x)+P^{-}_{2,2}(x) \\[1.5 ex]
\!\!&\!\!=\!\!&\!\!
-\frac{531441}{1435033600}\cdot x^{16}+\frac{177147}{20500480}\cdot x^{14}-\frac{59049}{394240}\cdot x^{12}+\frac{7003}{20160}\cdot x^{10}
\end{array}
$$
over $\left(0,\frac{\pi}{2}\right)$, i.e. it holds that
$$
f(x) > P(x)
$$
over $\left(0,\frac{\pi}{2}\right)$.

\medskip
\noindent
\textbf{IV} (The final part)
Based on the Sturm theorem, the following inequality
$$
P(x)>0
$$
is true over $\left(0,\frac{\pi}{2}\right)$.
The stated conclusion for the polynomial function $P$ is correct based on the following facts$:$\\
$1.$
We can conclude, by Sturm theorem, that the polynomial function $P(x)$ has only one zero over the concrete extended segment $[-0.1, 1.58]$
of the initial interval $\left(0,\frac{\pi}{2}\right)$.\\
$2.$
Facts
$P\left(-0.1\right) \neq 0$, $P\left(0\right) = 0$ and $P\left(1.58 \right) \neq 0$
are correct.\\
$3.$ The polynomial $P$ is positive at boundary point $\frac{\pi}{2}$
{\big (}$P\left(\frac{\pi}{2}\right) = 2.27261 \ldots > 0\,${\big )}.
$\!\!$\\[0.1 ex]
Therefore, the following inequality
$$
f(x)>0
$$
is true over  $\left(0,\frac{\pi}{2}\right)$. \hfill $\Box$}

%Fill author(s) affiliation(s), address(es) and emails here:
%example is shown
\noindent
\medskip

\vspace{1cc}
{\small
\noindent
\textbf{Bojan Banjac} \hfill(Received 08. 03. 2024.)\\
Computer Graphics Chair,  \hfill(Revised  19. 04. 2024.)\\
Faculty of Technical Sciences, University of Novi Sad, \\
Trg Dositeja Obradovi\' ca 16, Novi Sad, Serbia, \\
E-mail: {\it bojan.banjac@uns.ac.rs}

\vspace{0.1cc}
{\small
\noindent
\textbf{Branko Male\v {s}evi\' c} \\
Department of Applied Mathematics, \\
School of Electrical Engineering, University of Belgrade, \\
Bulevar kralja Aleksandra 73, Belgrade, Serbia, \\
E-mail: {\it branko.malesevic@etf.bg.ac.rs}

\vspace{0.1cc}
{\small
\noindent
\textbf{Milo\v s Mi\' covi\' c} \\
Department of Applied Mathematics, \\
School of Electrical Engineering, University of Belgrade, \\
Bulevar kralja Aleksandra 73, Belgrade, Serbia, \\
E-mail: {\it milos.micovic@etf.bg.ac.rs}

\vspace{0.1cc}
{\small
\noindent
\textbf{Bojana Mihailovi\' c} \\
Department of Applied Mathematics, \\
School of Electrical Engineering, University of Belgrade, \\
Bulevar kralja Aleksandra 73, Belgrade, Serbia, \\
E-mail: {\it mihailovicb@etf.bg.ac.rs}

\vspace{0.1cc}
{\small
\noindent
\textbf{Milica Savatovi\' c} \\
Department of Applied Mathematics, \\
School of Electrical Engineering, University of Belgrade, \\
Bulevar kralja Aleksandra 73, Belgrade, Serbia, \\
E-mail: {\it milica.makragic@etf.bg.ac.rs}

\newpage

\bigskip\medskip
\noindent
\centerline{\small \bf SUPPLEMENTARY MATERIAL}

\pagenumbering{gobble} 

\bigskip 

\medskip\noindent
Figure 1 illustrates the stratified family of functions 
$$
\hspace*{25 mm}
\varphi_p(x)
=
-\cos x \, + \, \left(\dfrac{\sin x}{x}\right)^{\!3} - \, \dfrac{1}{15} \, x^4 \, + \, p \, x^6
\hspace*{8 mm} {\Big (}\mbox{for $x \in \left(0,\dfrac{\pi}{2}\right)$}{\Big )}
$$
from Lemma 2.
Cases for all values of the parameter $p \!\in\! R^+$ are shown, highlighting those with constants obtained in Statement 1.

\medskip\noindent
\begin{figure}[hbt!]
\centering
\includegraphics[height=6.5cm]{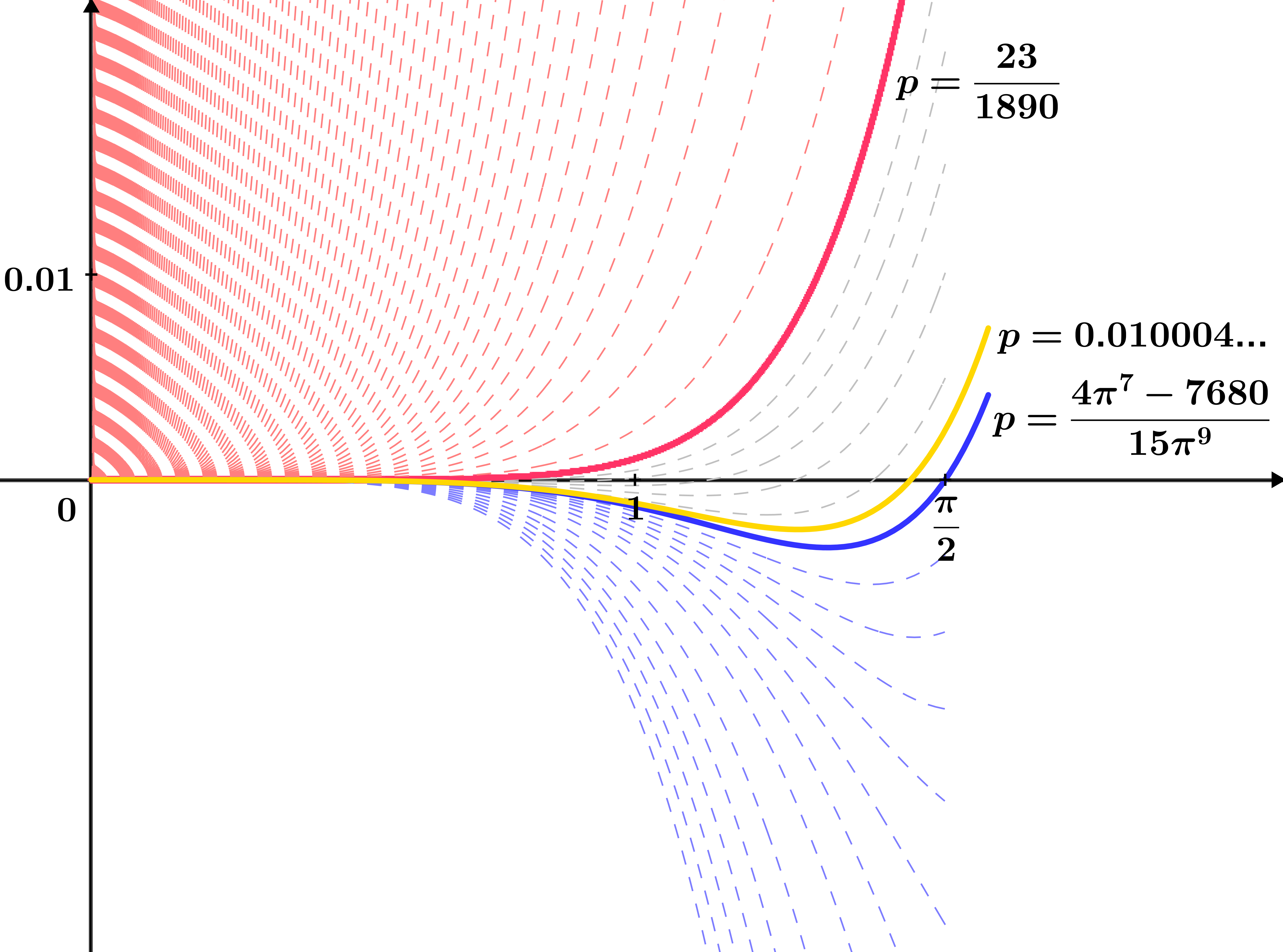}
\caption{Stratified family of functions from Lemma 2}
\end{figure}

\medskip\noindent
Figure 2 illustrates the stratified family of functions
$$
\hspace*{20.00 mm}
\varphi_p(x)
=
-\dfrac{\sin x}{x} + \dfrac{\cos x + 2}{3} - \dfrac{1}{180}x^4 + p \, x^6
\hspace*{10 mm} {\Big (}\mbox{for $x \in \left(0,\dfrac{\pi}{2}\right)$}{\Big )}
$$
from Lemma 4.
Cases for all values of the parameter $p \!\in\! R^+$ are shown, highlighting those with constants obtained in Statement 2.

\medskip\noindent
\begin{figure}[hbt!]
\centering
\includegraphics[height=6.5cm]{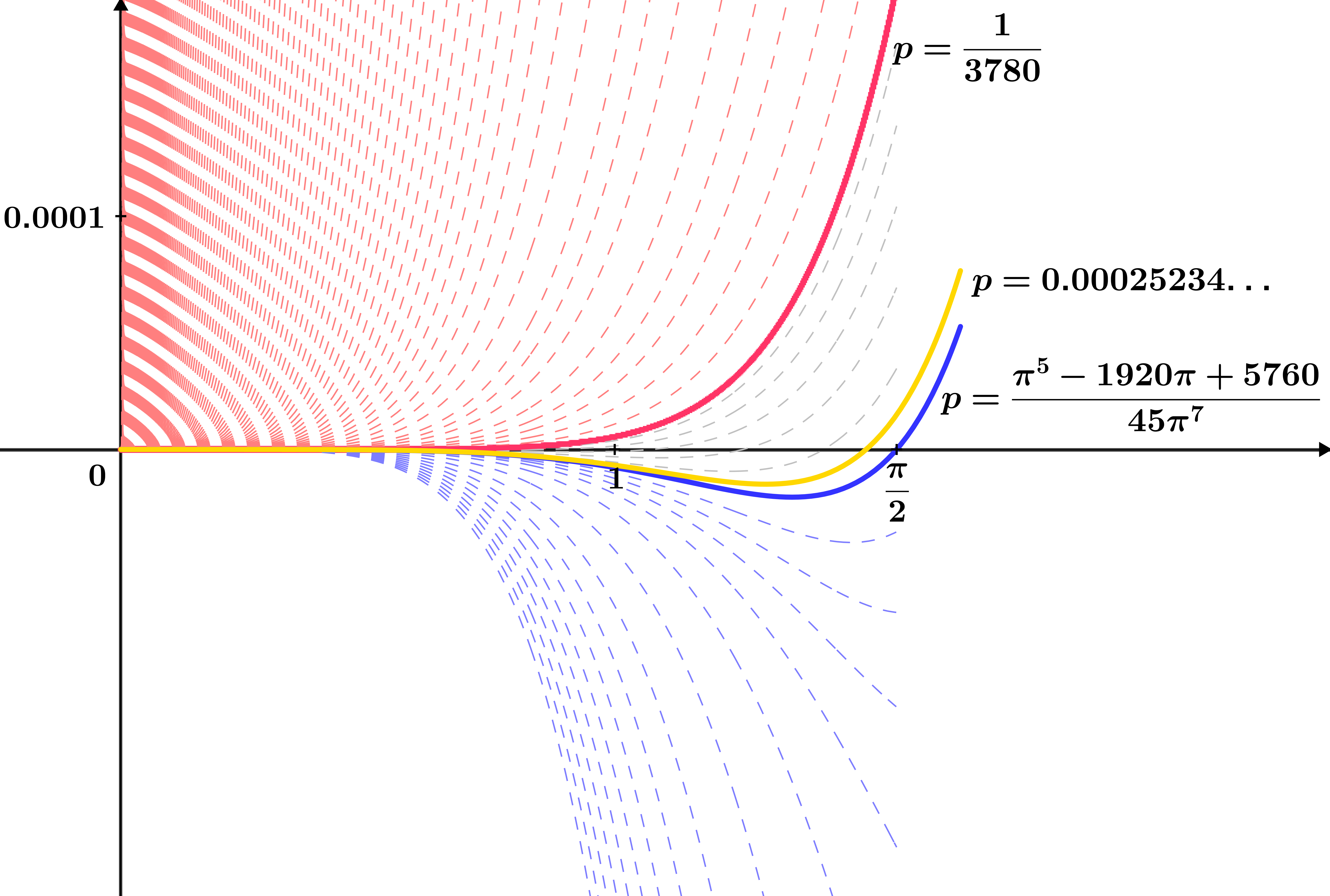}
\caption{Stratified family of functions from Lemma 4}
\end{figure}

\break

\medskip\noindent
Figure 3 illustrates the stratified family of functions
$$
\hspace*{5.00 mm}
\varphi_p(x)
=
2 \, \dfrac{\sin x}{x} + \dfrac{\tan x}{x} - 3 - \dfrac{3}{20} \, x^4 \! \dfrac{1}{\cos x} + p \, x^6 \dfrac{1}{\cos x}
\hspace*{2.25 mm} {\Big (}\mbox{for $x \in \left(0,\dfrac{\pi}{2}\right)$}{\Big )}
$$ 
from Lemma 6.
Cases for all values of the parameter $p \!\in\! R^+$ are shown, highlighting those with constants obtained in Statement 3.

\medskip\noindent
\begin{figure}[hbt!]
\centering
\includegraphics[height=6.5cm]{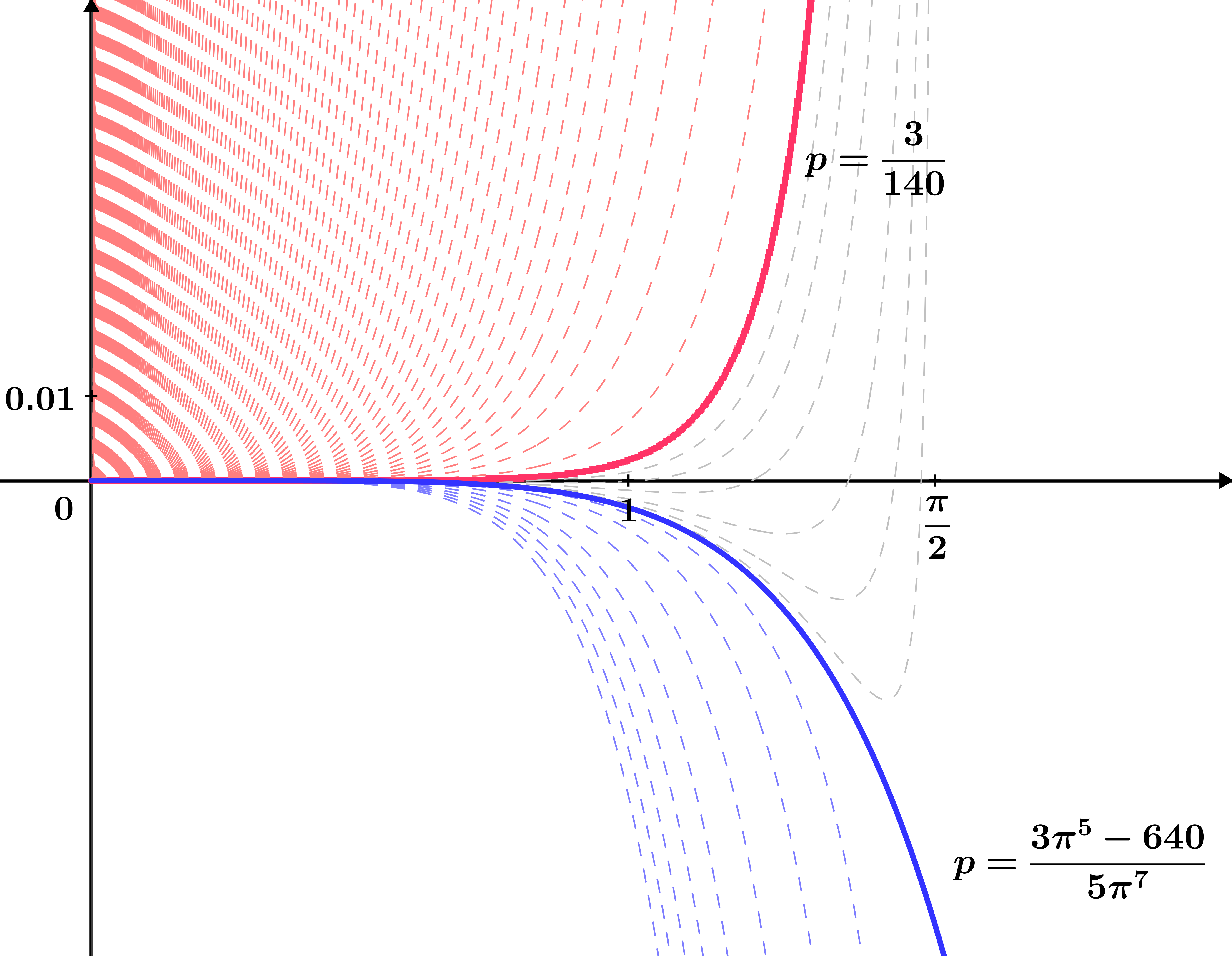}
\caption{Stratified family of functions from Lemma 6}
\end{figure}

\medskip\noindent
Figure 4 illustrates the stratified family of functions
$$
\hspace*{5.00 mm}
\varphi_p(x)
=
\left(\dfrac{\sin x}{x}\right)^2
+ \dfrac{\tan x}{x} - 2 - \dfrac{8}{45} \, x^4 \! \dfrac{1}{\cos x} + p \, x^6 \dfrac{1}{\cos x}
\hspace*{2.25 mm} {\Big (}\mbox{for $x \in \left(0,\dfrac{\pi}{2}\right)$}{\Big )}
$$ 
from Lemma 8.
Cases for all values of the parameter $p \!\in\! R^+$ are shown, highlighting those with constants obtained in Statement 4.

\medskip\noindent
\begin{figure}[hbt!]
\centering
\includegraphics[height=6.5cm]{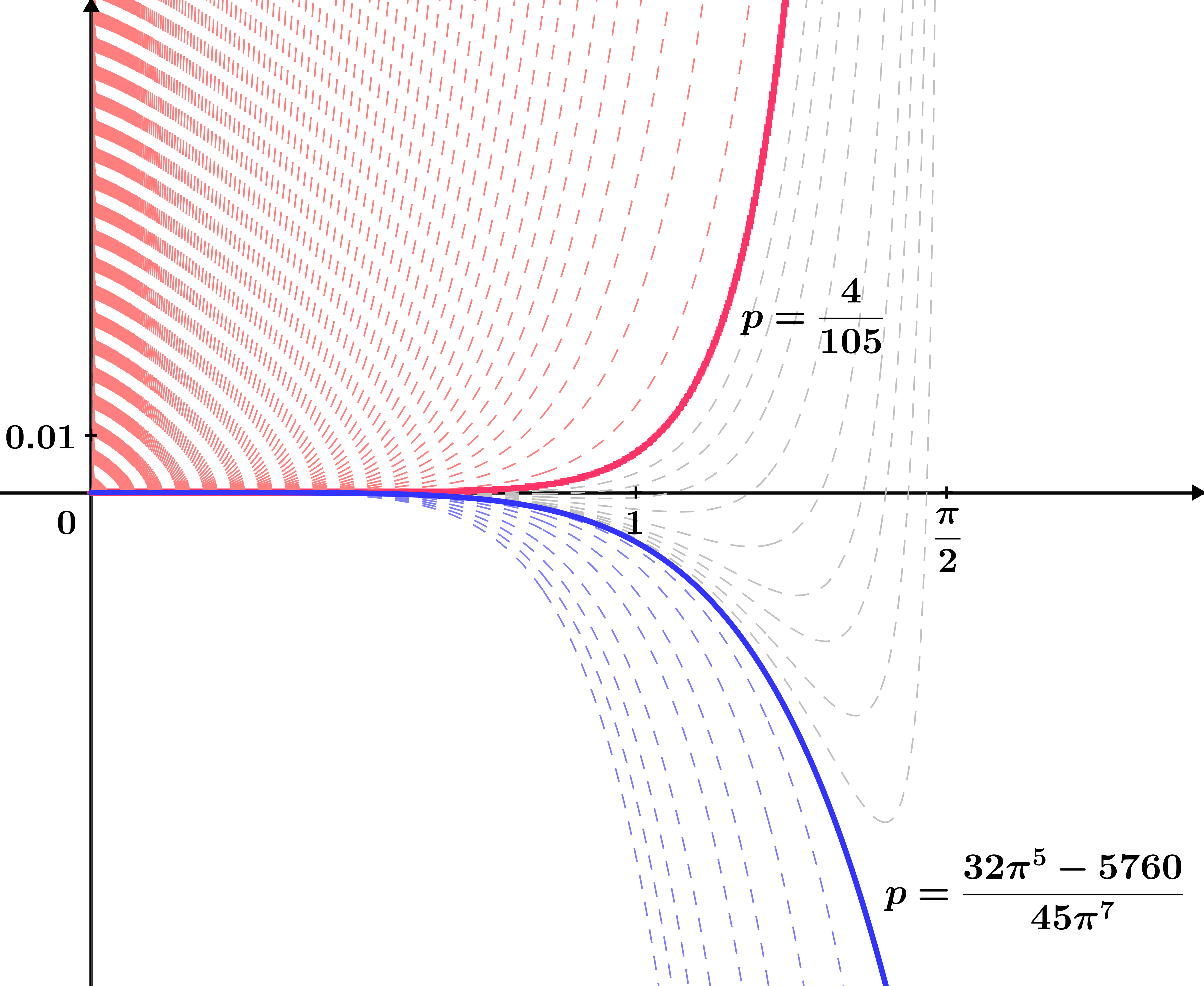}
\caption{Stratified family of functions from Lemma 8}
\end{figure}

\medskip\noindent
Figure 5 illustrates the stratified family of functions 
$$
\hspace*{35.0 mm}
\varphi_p(x)
=
\left(\dfrac{x}{\sin x}\right)^2
\, + \, \dfrac{x}{\tan x} - 2 - \, p \, x^4
\hspace*{8.80 mm} {\Big (}\mbox{for $x \in \left(0,\dfrac{\pi}{2}\right)$}{\Big )}
$$
from Lemma 10.
Cases for all values of the parameter $p \!\in\! R^+$ are shown, highlighting those with constants obtained in Statement 5.

\medskip\noindent
\begin{figure}[hbt!]
\centering
\includegraphics[height=6.5cm]{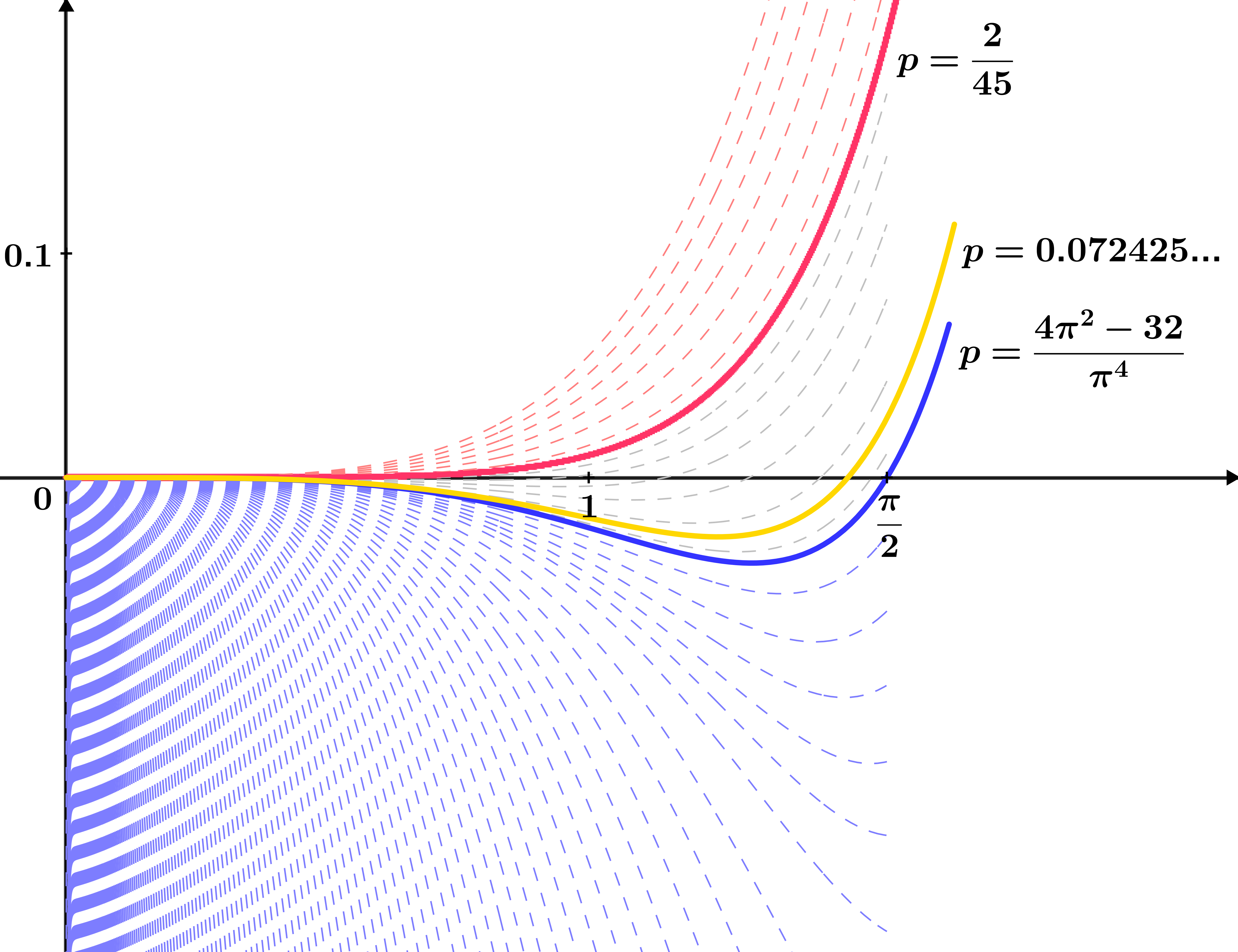}
\caption{Stratified family of functions from Lemma 10}
\end{figure}

\break

\medskip\noindent
Figure 6 illustrates the stratified family of functions 
$$
\hspace*{25.0 mm}
\varphi_p(x)
=
3 \, \dfrac{x}{\sin x}
\, + \, \cos x \, - \, 4 \, - \, \dfrac{1}{10} \, x^4 \, -  \, p \, x^6
\hspace*{8.80 mm} {\Big (}\mbox{for $x \in \left(0,\dfrac{\pi}{2}\right)$}{\Big )}
$$
from Lemma 11.
Cases for all values of the parameter $p \!\in\! R^+$ are shown, highlighting those with constants obtained in Statement 6.

\medskip\noindent
\begin{figure}[hbt!]
\centering
\includegraphics[height=6.5cm]{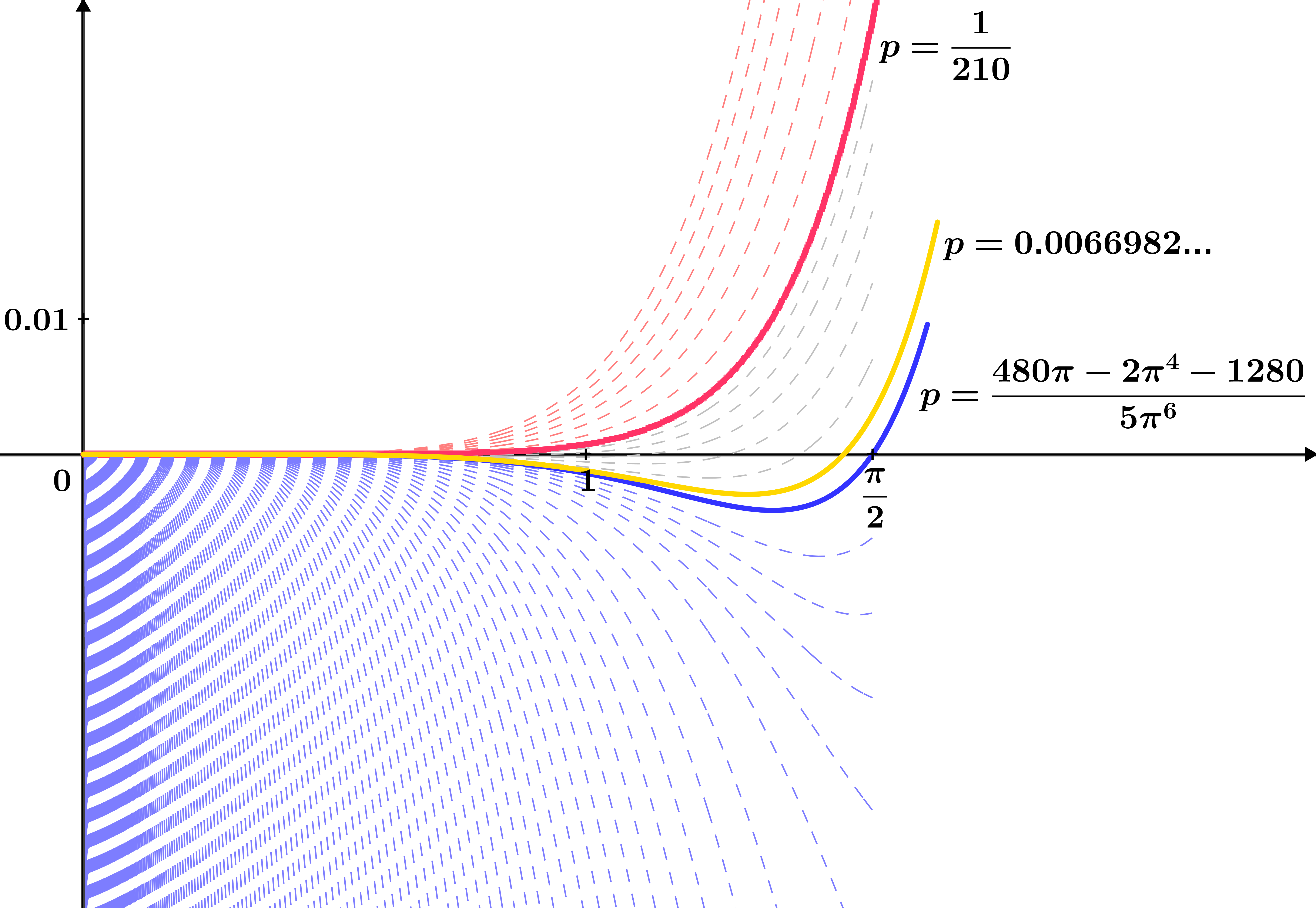}
\caption{Stratified family of functions from Lemma 11}
\end{figure}

\end{document}